\documentclass[a4paper, 12pt]{article}
\usepackage[russian]{babel}
\usepackage{amsmath}
\usepackage{amsfonts}

\newcommand{\mx}{\operatorname{\mathfrak M}}
\newcommand{\cmx}{\operatorname{\mathfrak C}}
\newcommand{\mn}{\operatorname{\mathfrak m}}
\newcommand{\cmn}{\operatorname{\mathfrak c}}
\newcommand{\e}{\operatorname{\mathfrak e}}
\newcommand{\s}{\operatorname{\mathsf s}}
\newcommand{\supp}{\operatorname{supp}}
\newcommand{\card}{\operatorname{card}}
\newcommand{\sgn}{\operatorname{sign}}
\newcommand{\supvrai}{\operatornamewithlimits{sup\,vrai}}
\newcommand{\N}{\mathbb N}
\newcommand{\Rho}{\mathrm P}
\newcommand{\Z}{\mathbb Z}
\newcommand{\R}{\mathbb R}
\newcommand{\Nu}{\mathcal N}
\newcommand{\D}{\mathcal D}

\begin{document}
\author{ С. Н. Кудрявцев }
\title{Верхняя оценка наилучшей точности восстановления производных
по значениям функций смешанной гладкости из классов Бесова}
\date{}\maketitle
\begin{abstract}
В работе установлена оценка сверху величины наилучшей точности
восстановления производных по значениям функций в заданном числе точек
для классов Бесова функций, удовлетворяющих
смешанным условиям Г\"ельдера. Эта оценка в некоторых случаях
является более сильной, чем соответствующая оценка,
полученная автором ранее в рассматриваемой задаче.
\end{abstract}

Ключевые слова: точность, восстановление, производная,  
значения функций, смешанная гладкость 
\bigskip

\centerline{Введение}

В работе при соответствующих условиях на $ \lambda \in \Z_+^d, d \in \N, $
установлена оценка сверху величины наилучшей
точности восстановления в пространстве $ L_q(I^d), 1 \le q \le \infty, $
производной $ \D^\lambda f $ по значениям в $ n $ точках функций $ f $
из классов Бесова $ \mathcal S_{p,\theta}^\alpha \mathcal B(I^d),
\alpha \in \R_+^d, 1 \le p < \infty, 1 \le \theta \le \infty. $

Настоящая работа продолжает исследования, проводившиеся автором в
[1] -- [4] для упомянутой задачи в отношении различных классов функций
конечной гладкости.

Отметим, что средства, использовавшиеся в [1] -- [3] для получения
как верхних, так и нижних оценок указанной выше величины для
рассматривавшихся в этих работах классов функций, неприменимы в
отношении классов, изучаемых в [4] и настоящей статье. Помимо этого,
схемы, применявшиеся в [1] -- [3] для вывода верхних оценок
точности приближений производных функций из исследуемых в [1] -- [3]
классов, оказались непригодными для изучаемых в [4] и настоящей работе
классов. Предложенные в [4] средства и схемы вывода оценок
позволили найти при определ\"енных условиях слабую асимптотику
отмеченной выше величины для классов Никольского и Бесова функций,
удовлетворяющих смешанным условиям Г\"ельдера.
Используя эти средства и применяя новые технические приемы вывода,
в настоящей работе получена в некоторых случаях более сильная оценка сверху
рассматриваемой здесь величины для таких классов функций.
Работа состоит из введения и двух параграфов, в первом из которых
излагаются предварительные сведения, а во втором -- основные резкльтаты.

В заключение в качестве иллюстрации привед\"ем полученный в \S 2 результат,
касающийся величины наилучшей точности
восстановления в $ L_q(I^d) $ производных $ \D^\lambda f $
по значениям в $ n $ точках функций $ f $ из класса Бесова
$ \mathcal S_{p, \theta}^\alpha \mathcal B(I^d), $
определяемой равенством
$$
\sigma_n(\D^\lambda,\mathcal S_{p,\theta}^\alpha \mathcal B(I^d),L_q(I^d))
= \inf_{A,\phi} \sup_{f \in \mathcal S_{p,\theta}^\alpha \mathcal B(I^d)}
\|\D^\lambda f - A \circ \phi(f)\|_{L_q(I^d)},
$$
где $ \inf $ бер\"ется по всем отображениям $ A: \R^n \mapsto
L_q(I^d) $ и отображениям $ \phi: C(I^d) \mapsto \R^n $ вида $
\phi(f) = (f(x^1),\ldots,f(x^n)), x^j \in I^d,j=1,\ldots.n. $
Для $ d \in \N, \alpha \in \R_+^d, \lambda \in \Z_+^d, 1 \le p < \infty,
1 \le \theta, q \le \infty $ установлено, что при $ \alpha_j -p^{-1} >0,
\alpha_j -\lambda_j -(p^{-1} -q^{-1})_+ >0, j=1,\ldots,d, $ справедливо
неравенство
$$
\sigma_n(\D^\lambda,\mathcal S_{p,\theta}^\alpha \mathcal B(I^d),L_q(I^d))
\le c_1 n^{-\mn} (\log n)^{(\mn +1 -1/\max(p,\theta))(\cmn -1)},
$$
где
\begin{multline*}
\mn = \min\{\alpha_j -\lambda_j -(p^{-1} -q^{-1})_+:
j=1,\ldots,d\}, \\
\cmn = \card\{j =1,\ldots,d: \alpha_j -\lambda_j -(p^{-1} -q^{-1})_+ = \mn \},
\end{multline*}

$ c_1 $ -- положительная константа, не зависящая от $ n. $
Отметим, что в аналогичной оценке из [4] слагаемое $ -1/\max(p,\theta) $
в показателе степени $ \log n $ отсутствует.
\bigskip

\centerline{\S 1. Предварительные сведения и вспомогательные
утверждения}
\bigskip

   1.1. В этом пункте вводятся обозначения, относящиеся к
функциональным классам и пространствам, рассматриваемым в
настоящей работе, а также приводятся некоторые факты, необходимые

в дальнейшем.

Для $ d \in \N $ через $ \Z_+^d $ обозначим множество
$$
\Z_+^d =\{\lambda =(\lambda_1,  \ldots, \lambda_d) \in \Z^d:
\lambda_j \ge0, j=1, \ldots, d\}.
$$
 Обозначим также при  $ d \in \N $ для $ l \in \Z_+^d $ через
$ \Z_+^d(l) $ множество
 $$
 \Z_+^d(l) =\{ \lambda  \in \Z_+^d: \lambda_j \le l_j, j=1, \ldots,  d\}.
 $$
Для $  d \in  \N, l \in \Z_+^d $ через $ \mathcal P^{d,l} $ будем
обозначать пространство вещественных  полиномов, состоящее из всех
функций $ f: \R^d \mapsto \R $ вида
$$
f(x) =\sum_{\lambda  \in \Z_+^d(l)}a_{\lambda}\cdot x^{\lambda},
x\in \R^d,
$$
где  $   a_{\lambda}   \in   \R,  x^{\lambda} =x_1^{\lambda_1}
\ldots x_d^{\lambda_d}, \lambda \in \Z_+^d(l). $

При $ d \in \N, l \in \Z_+^d  $   для области $ D\subset\R^d $
через $\mathcal   P^{d,l}(D) $  обозначим пространство функций $
f, $ определ\"енных в $  D, $   для каждой из  которых существует
полином $ g\in\mathcal P^{d,l} $
 такой, что  сужение $ g\mid_D = f.$

Для топологического пространства $ T $ через $ C(T) $
обозначим пространство непрерывных вещественных
функций, заданных на $ T. $ Если $ T $ --
компактное топологическое пространство, то в $ C(T) $ фиксируем норму,
определяемую равенством
$$
\|f\|_{C(T)} = \max_{x \in T} |f(x)| = \sup_{x \in T} |f(x)|,
f \in C(T).
$$

 Для измеримого по Лебегу множества $ D \subset \R^d$
при $ 1\le p\le\infty $ через  $  L_p(D),$  как обычно,
обозначается пространство  всех  вещественных измеримых на $ D $
функций $f,$ для которых определена норма
$$
\|f\|_{L_p(D)}= \begin{cases} (\int_D |f(x)|^p
dx)^{1/p}, 1 \le p<\infty;\\
\supvrai_{x \in D}|f(x)|, p=\infty. \end{cases}
$$

Введ\"ем ещ\"е следующие обозначения.

Для $ x,y \in \R^d $ положим $ xy =x \cdot y = (x_1 y_1, \ldots,
x_d y_d), $ а для $ x \in \R^d $ и $ A \subset \R^d $ определим
$$
x A = x \cdot A = \{xy: y \in A\}.
$$

Для $ x \in \R^d: x_j \ne 0, $ при $ j=1,\ldots,d,$ положим $
x^{-1} = (x_1^{-1},\ldots,x_d^{-1}). $

При $ d \in \N $ для $ x,y \in \R^d $ будем обозначать
$$
(x,y) = \sum_{j=1}^d x_j y_j.
$$

Обозначим через $ \R_+^d $ множество $ x \in \R^d, $ для которых $
x_j >0 $ при $ j=1,\ldots,d,$ и для $ a \in \R_+^d, x \in \R^d $
положим $ a^x = a_1^{x_1} \ldots a_d^{x_d}.$

При $ d \in \N $ определим множества
\begin{eqnarray*}
I^d &=& \{x \in \R^d: 0 < x_j < 1,j=1,\ldots,d\},\\
\overline I^d &=& \{x \in \R^d: 0 \le x_j \le 1,j=1,\ldots,d\},\\
B^d &=& \{x \in \R^d: -1 \le x_j \le 1,j=1,\ldots,d\}.
\end{eqnarray*}

Через $ \e $ будем обозначать вектор в $ \R^d, $ задаваемый
равенством $ \e =(1,\ldots,1). $

Далее, напомним, что для области $ D \subset \R^d $ и вектора $ h
\in \R^d $ через $ D_h $ обозначается множество
$$
D_h = \{x \in D: x +th \in D \forall t \in \overline I\}.
$$

Для функции $ f, $ заданной в области $ D \subset \R^d, $ и
вектора $ h \in \R^d $ определим в $ D_h $ е\"е разность $ \Delta_h
f $ с шагом $ h, $ полагая
$$
(\Delta_h f)(x) = f(x+h) -f(x), x \in D_h,
$$
а для $ l \in \N: l \ge 2, $ в $ D_{lh} $ определим $l$-ую
разность $ \Delta_h^l f $ функции $ f $ с шагом $ h $ равенством
$$
(\Delta_h^l f)(x) = (\Delta_h (\Delta_h^{l-1} f))(x), x \in
D_{lh},
$$
положим также $ \Delta_h^0 f = f. $

Как известно, справедливо равенство
$$
(\Delta_h^l f)(\cdot) = \sum_{k=0}^l C_l^k (-1)^{l-k} f(\cdot
+kh), C_l^k = l! /(k! (l-k)!).
$$

При $ d \in \N $ для $ j=1,\ldots,d$ через $ e_j $ будем
обозначать вектор $ e_j = (0,\ldots,0,1_j,0,\ldots,0).$

Лемма 1.1.1

   Пусть $ d \in \N, l \in \Z_+^d. $ Тогда
для любого $ \delta \in \R_+^d $ и $ x^0 \in \R^d $ для $ Q = x^0
+\delta I^d $ существует единственный линейный оператор $
\mathtt P_{\delta, x^0}^{d,l}: L_1(Q) \mapsto \mathcal P^{d,l}, $
обладающий следующими свойствами:

1) для $ f \in \mathcal P^{d,l} $ имеет место равенство
\begin{equation*}
\mathtt P_{\delta, x^0}^{d,l}(f \mid_Q) = f,
  \end{equation*}

2)
\begin{equation*}
\ker \mathtt P_{\delta,x^0}^{d,l} = \biggl\{\,f \in L_1(Q) : \int
\limits_{Q} f(x) g(x) \,dx =0\ \forall g \in \mathcal
P^{d,l}\,\biggr\},
\end{equation*}

прич\"ем, существуют константы $ c_1(d,l) >0 $ и $ c_2(d,l) >0 $
такие, что

   3) при $ 1 \le p \le \infty $ для $ f \in
L_p(Q) $ справедливо неравенство
  \begin{equation*}
\|\mathtt P_{\delta, x^0}^{d,l}f \|_{L_p(Q)} \le c_1 \|f\|_{L_p(Q)},
  \end{equation*}

4) при $ 1 \le p < \infty $ для $ f \in L_p(Q) $ выполняется
неравенство
   \begin{equation*}   \|f -\mathtt P_{\delta, x^0}^{d,l}f \|_{L_p(Q)} \le c_2 \sum_{j=1}^d
\delta_j^{-1/p}\biggl (\int_{\delta_j B^1} \int_{Q_{(l_j +1) \xi e_j}}
|\Delta_{\xi e_j}^{l_j +1} f(x)|^p dx d\xi\biggr)^{1/p}.
\end{equation*}

Доказательство леммы 1.1.1 приведено в [5].

При $ d \in \N $ для $ \lambda \in \Z_+^d $ через $ \D^\lambda $
будем обозначать оператор дифференцирования $ \D^\lambda =
\frac{\D^{|\lambda|}} {\D x_1^{\lambda_1} \ldots \D
x_d^{\lambda_d}}, $ где $ |\lambda| = \sum_{j=1}^d \lambda_j. $

В [3] содержится

    Лемма 1.1.2

Пусть $ d \in \N, l \in \Z_+^d, \lambda \in \Z_+^d(l), 1 \le p,q \le \infty,
\rho, \sigma \in \R_+^d. $ Тогда существует константа $ c_3(d,l,\lambda,\rho,
\sigma) >0 $ такая, что для любых измеримых по Лебегу множеств $ D,Q \subset
\R^d, $ для которых можно найти $ \delta \in \R_+^d $ и $ x^0 \in \R^d $ такие, что
$ D \subset (x^0 +\rho \delta B^d) $ и $ (x^0 +\sigma \delta I^d) \subset Q, $
для любого полинома $ f \in \mathcal P^{d,l} $ выполняется неравенство
    \begin{equation*} \tag{1.1.1}
\| \D^\lambda f\|_{L_q(D)} \le c_3 \delta^{-\lambda -p^{-1} \e +q^{-1} \e}
\|f\|_{L_p(q)}.
   \end{equation*}

Теперь определим классы функций, изучаемые в настоящей работе (см.
[6], [7]). Но прежде введ\"ем некоторые обозначения.

При $ d \in \N $ для $ x \in \R^d $ через $\s(x) $ обозначим
множество $ \s(x) = \{j =1,\ldots,d: x_j \ne 0\}, $ а для
множества $ J \subset \{1,\ldots,d\} $ через $ \chi_J $ обозначим
вектор из $ \R^d $ с компонентами
$$
(\chi_J)_j = \begin{cases} 1, &\text{для} j \in J; \\ 0,
&\text{для} j \in (\{1,\ldots,d\} \setminus J). \end{cases}
$$

При $ d \in \N $ для $ x \in \R^d $ и $ J = \{j_1,\ldots,j_k\}
\subset \N: 1 \le j_1 < j_2 < \ldots < j_k \le d, $ через $ x^J $
обозначим вектор $ x^J = (x_{j_1},\ldots,x_{j_k}) \in \R^k, $ а
для множества $ A \subset \R^d $ положим $ A^J = \{x^J: x \in A\}.$

Для области $ D \subset \R^d $ и векторов $ h \in \R^d $ и $ l \in
\Z_+^d $ через $ D_h^l $ обозначим множество
\begin{multline*}
D_h^l = (\ldots (D_{l_d h_d e_d})_{l_{d-1} h_{d-1} e_{d-1}}
\ldots)_{l_1 h_1 e_1}= \{ x \in D: x +tlh \in D \forall t \in \overline I^d\} \\
=\{ x \in D: (x +\sum_{j \in \s(l)} t_j l_j h_j e_j) \in D
\forall t^{\s(l)} \in (\overline I^d)^{\s(l)} \}.
\end{multline*}

Пусть $ d \in \N, D $ -- область в $ \R^d $ и $ 1 \le p \le
\infty. $ Тогда для $ f \in L_p(D), h \in \R^d $ и $ l \in \Z_+^d
$ определим в $ D_h^l $ смешанную разность  функции $ f $ порядка
$ l, $ соответствующую вектору $ h, $ равенством
\begin{multline*}
(\Delta_h^l f)(x) = ((\prod_{j=1}^d \Delta_{h_j e_j}^{l_j}) f)(x)
=((\prod_{j \in \s(l)} \Delta_{h_j e_j}^{l_j}) f)(x) \\
= \sum_{k \in \Z_+^d(l)} (-\e)^{l-k} C_l^k f(x+kh),\ x \in D_h^l,
\end{multline*}

где $ C_l^k = \prod_{j=1}^d C_{l_j}^{k_j}, k \in \Z_+^d(l). $

Имея в виду, что для $ f \in L_p(D), l \in \Z_+^d $ и векторов $
h,h^\prime \in \R^d: h^{\s(l)} = (h^\prime)^{\s(l)}, $ соблюдается
соотношение
$$
\| \Delta_h^l f\|_{L_p(D_h^l)} = \| \Delta_{h^\prime}^l
f\|_{L_p(D_{h^\prime}^l)},
$$
определим для функции $ f $ смешанный модуль непрерывности в $
L_p(D) $ порядка $ l $ равенством
$$
\Omega^l (f,t^{\s(l)})_{L_p(D)} = \supvrai_{ \{ h \in \R^d:
h^{\s(l)} \in t^{\s(l)} (B^d)^{\s(l)} \}} \| \Delta_h^l
f\|_{L_p(D_h^l)}, t^{\s(l)} \in (\R_+^d)^{\s(l)}.
$$

Пусть теперь $ \alpha \in \R_+^d, 1 \le p \le \infty $ и $ D $ --
область в $ \R^d. $ Тогда зададим вектор $ l = l(\alpha) \in \N^d,
$ полагая $ l_j = \min \{m \in \N: \alpha_j < m \}, j =1,\ldots,d,
$ и обозначим через $ S_p^\alpha H(D) (\mathcal S_p^\alpha
\mathcal H(D)) $ множество всех функций $ f \in L_p(D), $
обладающих тем свойством, что для любого непустого множества $ J
\subset \{1,\ldots,d\} $ выполняется неравенство
$$
\sup_{t^J \in
(\R_+^d)^J}(t^J)^{-\alpha^J}\Omega^{l\chi_J}(f,t^J)_{L_p(D)}
=\sup_{t^J \in (\R_+^d)^J}(\prod_{j \in J} t_j^{-\alpha_j})
\Omega^{l\chi_J}(f, t^{\s(l\chi_J)})_{L_p(D)}<\infty(\le1).
$$

Пусть $ \alpha,p,D $ и $ l =l(\alpha) $ -- те же, что и выше, и $
\theta \in \R: 1 \le \theta < \infty. $ Тогда обозначим через $
S_{p,\theta}^\alpha B(D) (\mathcal S_{p,\theta}^\alpha \mathcal
B(D)) $ множество всех функций $ f \in L_p(D), $ которые для
любого непустого множества $ J \subset \{1,\ldots,d\} $
удовлетворяют условию
\begin{multline*}
\biggl(\int_{(\R_+^d)^J} (t^J)^{-\e^J -\theta \alpha^J}
(\Omega^{l \chi_J}(f, t^J)_{L_p(D)})^\theta dt^J\biggr)^{1/\theta} \\
=\biggl(\int_{(\R_+^d)^J} (\prod_{j \in J} t_j^{-1 -\theta \alpha_j})
(\Omega^{l \chi_J}(f, t^{\s (l \chi_J)})_{L_p(D)})^\theta \prod_{j
\in J} dt_j\biggr)^{1/\theta} < \infty (\le 1).
\end{multline*}

При $ \theta = \infty $ положим $ S_{p,\infty}^\alpha B(D) =
S_p^\alpha H(D), \mathcal S_{p,\infty}^\alpha \mathcal B(D) =
\mathcal S_p^\alpha \mathcal H(D). $

Как известно (см., например, [4]), имеет место включение
$$
\mathcal S_{p, \theta}^\alpha \mathcal B(D) \subset c_4(\alpha)
\mathcal S_p^\alpha \mathcal H(D),
$$
где $ c_4(\alpha) = \prod_{j=1}^d 2^{1+\alpha_j}. $

В заключение этого пункта введ\"ем ещ\"е несколько обозначений.

Для банахова пространства $ X $ (над $ \R$) обозначим $ B(X) = \{x
\in X: \|x\|_X \le 1\}. $

Для банаховых пространств $ X,Y $ через $ \mathcal B(X,Y) $
обозначим банахово пространство, состоящее из непрерывных линейных
операторов $ T: X \mapsto Y, $ с нормой
$$
\|T\|_{\mathcal B(X,Y)} = \sup_{x \in B(X)} \|Tx\|_Y.
$$
Отметим, что если $ X=Y,$ то $ \mathcal B(X,Y) $ является
банаховой алгеброй.
\bigskip

1.2. В этом пункте содержатся сведения о кратных
рядах, которыми будем пользоваться в дальнейшем.

Для $ d \in \N, y \in \R^d $ положим
$$
\mn(y) = \min_{j=1,\ldots,d} y_j
$$
и для банахова пространства $ X, $ вектора $ x \in X $ и семейства
$ \{x_\kappa \in X, \kappa \in \Z_+^d\} $ будем писать $ x =
\lim_{ \mn(\kappa) \to \infty} x_\kappa, $ если для любого $
\epsilon >0 $ существует $ n_0 \in \N $ такое, что для любого $
\kappa \in \Z_+^d, $ для которого $ \mn(\kappa) > n_0, $
справедливо неравенство $ \|x -x_\kappa\|_X  < \epsilon. $

Пусть $ X $ -- банахово пространство, $ d \in \N $ и
$ \{ x_\kappa \in X, \kappa \in \Z_+^d\} $ -- семейство векторов.
Тогда под суммой ряда $ \sum_{\kappa \in \Z_+^d} x_\kappa $ будем
понимать вектор $ x \in X, $ для которого выполняется равенство $
x = \lim_{\mn(k) \to \infty} \sum_{\kappa \in \Z_+^d(k)} x_\kappa. $

При $ d \in \N $ через $ \Upsilon^d $ обозначим множество
$$
\Upsilon^d = \{ \upsilon \in \Z^d: \upsilon_j \in \{0,1\},
j=1,\ldots,d\}.
$$

Имеет место

   Лемма 1.2.1

Пусть $ X $ -- банахово пространство, а вектор $ x \in X $ и
семейство $ \{x_\kappa \in X, \kappa \in \Z_+^d\} $ таковы, что $
x = \lim_{ \mn(\kappa) \to \infty} x_\kappa, $ Тогда для семейства $
\{ \mathcal X_\kappa \in X, \kappa \in \Z_+^d \}, $ определяемого
равенством
$$
\mathcal X_\kappa = \sum_{\upsilon \in \Upsilon^d: \s(\upsilon)
\subset \s(\kappa)} (-\e)^\upsilon x_{\kappa -\upsilon}, \kappa
\in \Z_+^d,
$$
справедливо равенство
$$
x = \sum_{\kappa \in \Z_+^d} \mathcal X_\kappa.
$$

Лемма является следствием того, что при $ k \in \Z_+^d $
выполняется равенство
\begin{equation*}
\sum|{\kappa \in \Z_+^d(k)}
\mathcal X_\kappa = x_k \text{(см. [5])}.
\end{equation*}

Замечание.

Легко заметить, что для любого семейства чисел
$ \{x_\kappa \in \R: x_\kappa \ge 0, \kappa \in \Z_+^d\} , $ если ряд
$ \sum_{\kappa \in \Z_+^d} x_\kappa $ сходится, т.е. существует предел
$ \lim_{\mn(k) \to \infty} \sum_{\kappa \in \Z_+^d(k)} x_\kappa, $
что эквивалентно соотношению
$ \sup_{k \in \Z_+^d} \sum_{\kappa \in \Z_+^d(k)} x_\kappa < \infty, $

то для любой последовательности подмножеств
$ \{Z_n \subset \Z_+^d, n \in \Z_+\}, $
 таких, что $ \card Z_n < \infty,
Z_n \subset Z_{n+1},  n \in \Z_+, $
и $ \cup_{ n \in \Z_+} Z_n = \Z_+^d, $
справедливо равенство
$ \sum_{\kappa \in \Z_+^d} x_\kappa =
\lim_{ n \to \infty} \sum_{\kappa \in Z_n} x_\kappa =
\sup_{k \in \Z_+^d} \sum_{\kappa \in \Z_+^d(k)} x_\kappa. $
Отсюда несложно понять, что если для семейства векторов
$ \{x_\kappa \in X, \kappa \in \Z_+^d\}  $ банахова пространства $ X $
ряд $ \sum_{\kappa \in \Z_+^d} \| x_\kappa \|_X $ сходится,
то для любой последовательности подмножеств
$ \{Z_n \subset \Z_+^d, n \in \Z_+\}, $
 таких, что $ \card Z_n < \infty,
Z_n \subset Z_{n+1},  n \in \Z_+, $
и $ \cup_{ n \in \Z_+} Z_n = \Z_+^d, $
в $ X $ соблюдается равенство
$ \sum_{\kappa \in \Z_+^d} x_\kappa =
\lim_{ n \to \infty} \sum_{\kappa \in Z_n} x_\kappa. $

При $ d \in \N $ для $ x \in \R^d $ обозначим
\begin{eqnarray*}
\mx(x)  &=& \max_{j=1,\ldots,d} x_j, \\
\cmx(x) &=& \card \{j \in \{1,\ldots,d\}: x_j = \mx(x)\}, \\
\cmn(x) &=& \card \{j \in \{1,\ldots,d\}: x_j = \mn(x)\}.
\end{eqnarray*}

Лемма 1.2.2

Пусть $ d \in \N, \beta \in \R_+^d, \alpha \in \R^d $ и $
\mx(\beta^{-1} \alpha) >0. $ Тогда существуют константы $
c_1(d,\alpha,\beta)  >0 $ и $ c_2(d,\alpha,\beta) >0 $ такие, что
для $ r \in \N $ соблюдается неравенство
\begin{equation*} \tag{1.2.1}
c_1 2^{ \mx(\beta^{-1} \alpha) r} r^{ \cmx(\beta^{-1} \alpha) -1}
\le \sum_{ \kappa \in \Z_+^d: (\kappa, \beta) \le r} 2^{(\kappa,
\alpha)} \le c_2 2^{ \mx(\beta^{-1} \alpha) r} r^{ \cmx(\beta^{-1}
\alpha) -1}.
\end{equation*}

Лемма 1.2.3

Пусть $ d \in \N, \alpha, \beta \in \R_+^d. $ Тогда существуют
константы $ c_3(d,\alpha,\beta) >0 $ и $ c_4(d,\alpha,\beta) >0 $
такие, что при $ r \in \N $ справедливо неравенство
\begin{equation*} \tag{1.2.2}
c_3 2^{-\mn(\beta^{-1} \alpha) r} r^{\cmn(\beta^{-1} \alpha) -1}
\le \sum_{\kappa \in \Z_+^d: (\kappa, \beta) > r} 2^{-(\kappa,
\alpha)} \le c_4 2^{-\mn(\beta^{-1} \alpha) r} r^{\cmn(\beta^{-1}
\alpha) -1}.
\end{equation*}

Доказательство лемм 1.2.2 и 1.2.3 приведено в [8].
\bigskip

1.3. В этом пункте привед\"ем некоторые сведения о рассматриваемых в работе
функциональных пространствах и действующих в них операторах, которые
используются в следующем параграфе.

Как показано в  [3], [4], справедлива

   Лемма 1.3.1

    Пусть $ d \in \N, 1 \le p < \infty. $ Тогда

1) при  $ j=1,\ldots,d$ для любого непрерывного линейного
оператора $ T: L_p(I) \mapsto L_p(I) $ существует единственный
непрерывный линейный оператор $ \mathcal T^j: L_p(I^d) \mapsto
L_p(I^d), $ для которого для любой функции $ f \in L_p(I^d) $
почти для всех $ (x_1,\ldots,x_{j-1},x_{j+1}, \ldots,x_d) \in
I^{d-1} $ в $ L_p(I) $ выполняется равенство
\begin{multline*} \tag{1.3.1}
(\mathcal
T^j f)(x_1,\ldots,x_{j-1},\cdot,x_{j+1},\ldots,x_d)\\
=(T(f(x_1,\ldots,x_{j-1},\cdot,x_{j+1},\ldots,x_d)))(\cdot),
\end{multline*}

2) при этом, для каждого $ j=1,\ldots,d $ отображение $ V_j^{L_p},
$ которое каждому оператору $ T \in \mathcal B(L_p(I), L_p(I)) $
ставит в соответствие оператор $ V_j^{L_p}(T) = \mathcal T^j \in
\mathcal B(L_p(I^d), L_p(I^d)), $ удовлетворяющий (1.3.1),
является непрерывным гомоморфизмом банаховой алгебры $ \mathcal
B(L_p(I), L_p(I)) $ в банахову алгебру $ \mathcal B(L_p(I^d),
L_p(I^d)), $

3) прич\"ем, для любых операторов $ S,T \in \mathcal B(L_p(I),
L_p(I)) $ при любых $ i,j =1,\ldots,d: i \ne j, $ соблюдается
равенство
\begin{equation*}
(V_i^{L_p}(S) V_j^{L_p}(T))f = (V_j^{L_p}(T) V_i^{L_p}(S))f, f \in
L_p(I^d).
\end{equation*}

В [4] также установлена

     Лемма 1.3.2

    Пусть $ d \in \N. $ Тогда

1) при  $ j=1,\ldots,d$ для любого непрерывного линейного оператора
$ T: C(\overline I)
\mapsto C(\overline I) $ существует единственный непрерывный линейный
оператор $ \mathcal T^j: C(\overline I^d) \mapsto C(\overline I^d), $
обладающий тем свойством, что для любой функции
$ f \in C(\overline I^d) $ для всех $ (x_1,\ldots,x_{j-1},x_{j+1},
\ldots,x_d) \in \overline I^{d-1} $ при любом $ x_j \in \overline I $
имеет место равенство
\begin{multline*} \tag{1.3.2}
(\mathcal T^j f)(x_1,\ldots, x_{j-1},x_j,x_{j+1},
\ldots,x_d) \\
= (T(f(x_1,\ldots,x_{j-1},\cdot,x_{j+1}, \ldots,x_d)))(x_j),
\end{multline*}

2) при этом, для каждого $ j=1,\ldots,d $ отображение $ V_j^C, $ которое каждому
оператору $ T \in \mathcal B(C(\overline I), C(\overline I)) $ ставит в
соответствие оператор $ V_j^C(T) = \mathcal T^j \in \mathcal B(C(\overline I^d),
C(\overline I^d)), $
удовлетворяющий (1.3.2), является непрерывным гомоморфизмом банаховой алгебры
$ \mathcal B(C(\overline I), C(\overline I)) $ в банахову алгебру $ \mathcal
B(C(\overline I^d), C(\overline I^d)), $

3) прич\"ем, для любых операторов $ S,T \in \mathcal B(C(\overline I),
C(\overline I)) $
при любых $ i,j =1,\ldots,d: i \ne j, $ справедливо равенство
\begin{equation*} \tag{1.3.3}
V_i^C(S) V_j^C(T) = V_j^C(T) V_i^C(S).
\end{equation*}

Замечание.

Если при $ d \in \N, 1 \le p \le q < \infty, $ оператор
$ T \in \mathcal B(L_p(I), L_p(I)) \cap \mathcal B(L_q(I), L_q(I))
( T \in \mathcal B(L_p(I), L_p(I)) \cap \mathcal B(C(\overline I),
C(\overline I))), $ то при $ j=1,\ldots,d $ для $ f \in L_q(I^d)
(f \in C(\overline I^d)) $
справедливо равенство $ (V_j^{L_p} T)f = (V_j^{L_q} T)f
((V_j^{L_p} T)f = (V_j^C T)f). $
Поэтому символы $ L_p, L_q, C $ в качестве индексов у $ V_j $ можно опускать.

Привед\"ем ещ\"е одно вспомогательное утверждение.

  Лемма 1.3.3

  Пусть $ d \in \N, l \in \N^d, 1 \le p < \infty. $ Тогда существует константа
$ c_1(d,l) >0 $ такая, что для любых $ x^0 \in \R^d, \delta \in \R_+^d $
таких, что $ Q = (x^0 +\delta I^d) \subset I^d, $
для любой функции $ f \in L_p(I^d), $ для любого $ \xi \in \R^d, $ для любых
множеств $ J,J^\prime \subset \{1,\ldots,d\} $
имеет место неравенство
\begin{multline*} \tag{1.3.4}
\| \Delta_\xi^{l \chi_J} ((\prod_{j \in J^\prime}
V_j(E- \mathtt P_{\delta_j, x_j^0}^{1,l_j -1}))
f)\|_{L_p(Q_\xi^{l \chi_J})} \\
\le c_1 \biggl(\prod_{j \in J^\prime \setminus J}
\delta_j^{-1/p}\biggr) \biggl(\int_{(\delta B^d)^{J^\prime \setminus J}}
\int_{ Q_\xi^{l \chi_{J \cup J^\prime}}}
|(\Delta_\xi^{l \chi_{J \cup J^\prime}} f)(x)|^p dx
d\xi^{J^\prime \setminus J}\biggr)^{1/p},
\end{multline*}
где $ E $ -- тождественный оператор.

Доказательство леммы 1.3.3 дословно повторяет доказательство аналогичного
утверждения в [4]. В [4] также доказаны теорема 1.3.4 и предложение 1.3.5.

Теорема 1.3.4

Пусть $ d \in \N, \alpha \in \R_+^d, 1 \le p < \infty, 1 \le q \le
\infty $ и $ \lambda \in \Z_+^d $ удовлетворяют условию
\begin{equation*} \tag{1.3.5}
\alpha -\lambda -(p^{-1} -q^{-1})_+ \e >0.
\end{equation*}
Тогда существуют константы $ c_{2}(d,\alpha) >0$
и $ c_{3}(d,\alpha,p,q,\lambda) >0 $ такие, что
при любых $ \delta \in \R_+^d, x^0 \in \R^d $ для
$ Q = x^0 +\delta I^d $ для $ f \in S_p^\alpha H(Q) $ при $ l =
l(\alpha) $ имеет место неравенство
\begin{multline*} \tag{1.3.6}
\| \D^\lambda f \|_{L_q(Q)} \le c_{3} \delta^{-\lambda -p^{-1} \e
+q^{-1} \e} \biggl(\|f\|_{L_p(Q)} + \sum_{ J \subset \Nu_{1,d}^1:
J \ne \emptyset} (\prod_{j \in J} \delta_j^{\lambda_j +(p^{-1}
-q^{-1})_+}) \\ \times \int_{ (c_{2} \delta I^d)^J} (\prod_{j \in
J} t_j^{-\lambda_j -p^{-1} -(p^{-1} -q^{-1})_+ -1})
\biggl(\int_{(t B^d)^J} \int_{ Q_\xi^{l \chi_J}} |\Delta_\xi^{l
\chi_J} f(x)|^p dx d\xi^J\biggr)^{1/p} dt^J\biggr).
\end{multline*}

Предложение 1.3.5

Пусть $ d \in \N, \alpha \in \R_+^d, 1 \le p < \infty $
таковы, что
\begin{equation*} \tag{1.3.7}
\alpha -p^{-1} \e >0.
\end{equation*}
Тогда для любой функции $ f \in S_p^\alpha H(I^d) $ существует
функция $ F \in C(\overline I^d) $ такая, что почти для всех $ x
\in I^d $  верно равенство $ f(x) = F(x), $ т.е. $ S_p^\alpha
H(I^d) \subset C(\overline I^d). $
 \bigskip

\centerline{\S 2. Оценка сверху наилучшей точности восстановления
в $ L_q(I^d) $ производной $ \D^\lambda f$}
\centerline{по значениям в $n$ точках функций $ f$  из $ \mathcal S_{p,\theta}^\alpha \mathcal B(I^d) $}

2.1. В этом пункте определим объекты, используемые
для построения подходящих приближений, и приведем ещ\"е факты,
на которые опирается вывод интересующей нас оценки.
Введ\"ем в рассмотрение систему разбиений единицы на кубе $ I^d, $
используемую для построения средств приближения.
Для этого обозначим через $ \psi^{1,0} $ характеристическую
функцию интервала $ I, $ т.е. функцию определяемую равенством
$$
\psi^{1,0}(x) = \begin{cases} 1, & \text{для} x \in I; \\
0, & \text{для} x \in \R \setminus I.
\end{cases}
$$
При $ m \in \N $ положим
$$
\psi^{1,m}(x) = \int_I \psi^{1, m-1}(x-y) dy,
$$
а для $ d \in \N, m \in \Z_+^d $ определим
$$
\psi^{d,m}(x) = \prod_{j=1}^d \psi^{1,m_j}(x_j),
x = (x_1,\ldots,x_d) \in \R^d.
$$

Для $ d \in \N, x,y \in \R^d $ будем писать $ x \le y (x < y), $
если для каждого $ j=1,\ldots,d $ выполняется неравенство $ x_j \le y_j
(x_j < y_j). $

Для $ d \in \N, m,n \in \Z^d: m \le n, $ обозначим
$$
\Nu_{m,n}^d = \{ \nu \in \Z^d: m \le \nu \le n \} =
\prod_{j=1}^d \Nu_{m_j,n_j}^1.
$$

Опираясь на определения, используя индукцию, нетрудно проверить следующие
свойства функций $ \psi^{d,m}, d \in \N, m \in \Z_+^d, $

1) При $ d \in \N, m \in \Z_+^d $
$$
\sgn \psi^{d,m}(x) = \begin{cases} 1, \text{для} x \in ((m+\e)I^d); \\
0, \text{для} x \in \R^d \setminus ((m+\e) I^d),
\end{cases}
$$

2) при $ d \in \N, m \in \Z_+^d $ для каждого $ \lambda \in \Z_+^d(m) $
обобщ\"енная производная $ \D^\lambda \psi^{d,m} \in L_\infty(\R^d), $

3) при $ d \in \N, m \in \Z_+^d $ почти для всех $ x \in \R^d $ справедливо равенство
$$
\sum_{\nu \in \Z^d} \psi^{d,m}(x -\nu) =1,
$$

4) при $ m \in \Z_+ $ почти для всех $ x \in \R $ имеет место равенство
\begin{equation*}
\psi^{1,m}(x) = \sum_{\mu \in \Nu_{0, m+1}^1} a_{\mu}^m \psi^{1,m}(2x -\mu),
\end{equation*}
где $ a_\mu^m = 2^{-m} C_{m+1}^\mu, $

5) при $ m \in \Z_+ $ система функций $ \{ \psi^{1,m}(x -\nu), \nu \in
\Nu_{-m, 0}^1\} $ -- линейно независима на $ I. $

При $ d \in \N $ для $ t \in \R^d $ через $ 2^t $
будем обозначать вектор $ 2^t = (2^{t_1}, \ldots, 2^{t_d}). $

Для $ d \in \N, m,\kappa \in \Z_+^d, \nu \in \Z^d $ обозначим
$$
g_{\kappa, \nu}^{d,m}(x) = \psi^{d,m}(2^\kappa x -\nu) = \prod_{j=1}^d
\psi^{1,m_j}( 2^{\kappa_j} x_j -\nu_j), x \in \R^d,
$$
$$
Q_{\kappa, \nu}^d = 2^{-\kappa} \nu +2^{-\kappa} I^d.
$$

Из первого среди приведенных выше свойств функций
$ \psi^{d,m} $ следует, что при $ d \in \N, m,\kappa \in \Z_+^d, \nu \in
\Z^d $ носитель $ \supp g_{\kappa,\nu}^{d,m}
= 2^{-\kappa} \nu +2^{-\kappa} (m+\e)
\overline I^d. $

Отметим некоторые полезные для нас свойства носителей функций
$ g_{\kappa,\nu}^{d,m}. $

Для $ d \in \N, m,\kappa \in \Z_+^d $ справедливо
следующее:
\begin{equation*} \tag{2.1.1}
 \{\nu \in \Z^d: \supp g_{\kappa,\nu}^{d,m} \cap I^d \ne
\emptyset \} = \Nu_{-m, 2^\kappa -\e}^d,
\end{equation*}

2) для каждого $ \nu^\prime \in \Z^d $ имеет место равенство
\begin{equation*} \tag{2.1.2}
\{ \nu \in \Z^d: Q_{\kappa, \nu^\prime}^d \cap \supp
g_{\kappa,\nu}^{d,m} \ne \emptyset \} = \nu^\prime +\Nu_{-m,0}^d.
\end{equation*}

Из (2.1.2) вытекает, что
\begin{equation*} \tag{2.1.3}
\card \{ \nu \in \Nu_{-m, 2^\kappa -\e}^d: \supp g_{\kappa,
\nu}^{d,m} \cap Q_{\kappa, \nu^\prime}^d \ne \emptyset \} \le
c_1(d,m),
\end{equation*}
$$
d \in \N, m,\kappa \in \Z_+^d, \nu^\prime \in \Nu_{0, 2^\kappa
-\e}^d.
$$

Из свойства 3) функций $ \psi^{d,m} $ и равенства (2.1.1) вытекает, что при
$ d \in \N, m,\kappa \in \Z_+^d $ почти для всех $ x \in I^d $
соблюдается равенство
\begin{equation*}
\sum_{ \nu \in \Nu_{-m, 2^\kappa -\e}^d} g_{\kappa, \nu}^{d,m}(x) =1.
\end{equation*}

Имея в виду свойство 2) функций $ \psi^{d,m}, $ отметим, что при
$ d \in \N, m,\kappa \in \Z_+^d, \nu \in \Z^d,
\lambda \in \Z_+^d(m) $ выполняется равенство
\begin{equation*} \tag{2.1.4}
\| \D^\lambda g_{\kappa, \nu}^{d,m} \|_{L_\infty (\R^d)} =
c_2(d,m,\lambda) 2^{(\kappa, \lambda)}.
\end{equation*}

При доказательстве леммы 2.2.1 понадобится

    Лемма 2.1.1

Пусть $ d \in \N, \lambda \in \Z_+^d, D $ --- область в $ \R^d $ и
$ f $ --- бесконечно дифференцируемая в области $ D $ функция,
а $ g \in  L_1(D), $ прич\"ем,
для каждого $ \mu \in \Z_+^d(\lambda) $ обобщ\"енная  производная  $
\D^\mu g \in L_1(D). $ Тогда в пространстве обобщ\"енных функций в
области $ D $ имеет место равенство
\begin{equation*} \tag{2.1.5}
\D^\lambda (fg) = \sum_{ \mu  \in  \Z_+^d(\lambda)} C_\lambda^\mu
\D^{\lambda -\mu}f \D^\mu g.
\end{equation*}

Теперь напомним некоторые сведения, касающиеся интерполяционных полиномов.

При $ l \in \Z_+ $ в интервале $ I $ фиксируем систему из $ (l+1)
$ различных точек $ \{ \xi_{1,0}^{1,l,\lambda} \in I, \lambda = 0,
\ldots, l \} $ и построим систему полиномов $ \{
\pi_{1,0}^{1,l,\lambda} \in \mathcal P^{1,l}, \lambda =0, \ldots,
l \}, $ обладающую тем свойством, что при $ \lambda, \mu =0,
\ldots, l $ выполняются равенства
$$
\pi_{1,0}^{1,l,\lambda}(\xi_{1,0}^{1,l,\mu}) = \begin{cases} 1,
\text{при} \lambda = \mu; \\ 0, \text{при} \lambda \ne \mu.
\end{cases}
$$

При $ l \in \Z_+, \delta \in \R_+, x^0 \in \R $ определим систему
точек $ \{ \xi_{\delta, x^0}^{1,l,\lambda} \in x^0 +\delta I,
\lambda =0,\ldots,l\} $ и систему полиномов $ \{ \pi_{\delta,
x^0}^{1,l,\lambda} \in \mathcal P^{1,l}, \lambda =0,\ldots,l\}, $
полагая
$$
\xi_{\delta, x^0}^{1,l,\lambda} =
x^0 +\delta \xi_{1,0}^{1,l,\lambda}, \\
\pi_{\delta, x^0}^{1,l,\lambda}(x)  =
\pi_{1,0}^{1,l,\lambda}(\delta^{-1} (x -x^0)).
$$
Ясно, что при $ \lambda, \mu =0, \ldots,l $ выполняются равенства
\begin{equation*} \tag{2.1.6}
\pi_{\delta, x^0}^{1,l,\lambda}(\xi_{\delta, x^0}^{1,l,\mu}) =
\pi_{1,0}^{1,l,\lambda}(\xi_{1,0}^{1,l,\mu}) = \begin{cases} 1,
\text{при} \lambda = \mu; \\ 0, \text{при} \lambda \ne \mu.
\end{cases}
\end{equation*}

При $ d \in \N, l \in \Z_+^d, \delta \in \R_+^d, x^0 \in \R^d $
построим систему точек $ \{ \xi_{\delta, x^0}^{d,l,\lambda} \in
x^0 +\delta I^d, \lambda \in \Z_+^d(l) \} $ и систему полиномов $
\{ \pi_{\delta, x^0}^{d,l,\lambda} \in \mathcal P^{d,l}, \lambda
\in \Z_+^d(l) \}, $ задавая
$$
( \xi_{\delta, x^0}^{d,l,\lambda})_j =
 \xi_{\delta_j, x_j^0}^{1,l_j,\lambda_j}, j=1,\ldots,d,\\
 \pi_{\delta, x^0}^{d,l,\lambda}(x) =
\prod_{j=1}^d
 \pi_{\delta_j, x_j^0}^{1,l_j,\lambda_j}(x_j).
$$

Из (2.1.6) следует, что при $ \lambda, \mu \in \Z_+^d(l) $
справедливы равенства
\begin{equation*} \tag{2.1.7}
\pi_{\delta, x^0}^{d,l,\lambda}(\xi_{\delta, x^0}^{d,l,\mu}) =
\prod_{j=1}^d \pi_{\delta_j,
x_j^0}^{1,l_j,\lambda_j}(\xi_{\delta_j, x_j^0}^{1,l_j,\mu_j}) =
 \begin{cases} 1,
\text{при} \lambda = \mu; \\ 0, \text{при} \lambda \ne \mu.
\end{cases}
\end{equation*}

Понятно, что
\begin{multline*}
\pi_{\delta, x^0}^{d,l,\lambda}(x) = \prod_{j=1}^d
 \pi_{1, 0}^{1,l_j,\lambda_j}(\delta_j^{-1} (x_j -x_j^0)) =
\pi_{\e, 0}^{d,l,\lambda}(\delta^{-1} (x -x^0)), \\
\xi_{\delta, x^0}^{d,l,\lambda} = x^0 +\delta \xi_{\e,
0}^{d,l,\lambda}, \lambda \in \Z_+^d(l).
\end{multline*}

При $ d \in \N, l \in \Z_+^d, \delta \in \R_+^d, x^0 \in \R^d $
обозначим через $ A_{\delta, x^0}^{d,l}: \R^{ (l+\e)^{\e} }
\mapsto \mathcal P^{d,l} $ линейный оператор, который каждому
набору чисел $ t = \{ t_\lambda \in \R, \lambda \in \Z_+^d(l) \} $
ставит в соответствие полином
$$
A_{\delta, x^0}^{d,l} t = \sum_{ \lambda \in \Z_+^d(l)} t_\lambda
\pi_{\delta, x^0}^{d,l,\lambda} \in \mathcal P^{d,l},
$$
а через $ \phi_{\delta, x^0}^{d,l} $ обозначим непрерывное
линейное отображение пространства $ C(x^0 +\delta \overline I^d) $
в $ \R^{(l +\e)^{\e}}, $ которое каждой функции $ f \in C(x^0
+\delta \overline I^d) $ сопоставляет набор е\"е значений $ \{
f(\xi_{\delta, x^0}^{d,l,\lambda}), \lambda \in \Z_+^d(l) \}, $ и
определим линейный непрерывный оператор $ \mathcal P_{\delta,
x^0}^{d,l}: C(x^0 +\delta \overline I^d) \mapsto \mathcal P^{d,l}
\cap C(x^0 +\delta \overline I^d) $ равенством
$$
 \mathcal P_{\delta, x^0}^{d,l} =
A_{\delta, x^0}^{d,l} \circ \phi_{\delta, x^0}^{d,l}.
$$

Как видно из (2.1.7), при $ \lambda \in \Z_+^d(l) $ имеет место
равенство
\begin{equation*}
(\mathcal P_{\delta, x^0}^{d,l} f) (\xi_{\delta,
x^0}^{d,l,\lambda}) = f(\xi_{\delta, x^0}^{d,l,\lambda}).
\end{equation*}

Заметим, что вследствие (2.1.7) система полиномов $ \{
\pi_{\delta, x^0}^{d,l,\lambda} \in \mathcal P^{d,l}, \lambda \in
\Z_+^d(l) \} $ -- линейно независима, и число элементов этой
системы $ \card \Z_+^d(l) = \dim \mathcal P^{d,l}, $ и,
следовательно, эта система является базисом в $ \mathcal P^{d,l}.
$ Поэтому для любого полинома $ f \in \mathcal P^{d,l} $
существует набор чисел $ \{ t_\lambda \in \R, \lambda \in
\Z_+^d(l) \} $ такой, что
$$
f = \sum_{ \lambda \in \Z_+^d(l)} t_\lambda \pi_{\delta,
x^0}^{d,l,\lambda}.
$$
Отсюда в виду (2.1.7) получаем
$$
t_\lambda = f(\xi_{\delta, x^0}^{d,l,\lambda}), \lambda \in
\Z_+^d(l),
$$
т.е.
\begin{equation*} \tag{2.1.8}
f = \mathcal P_{\delta, x^0}^{d,l} f, f \in \mathcal P^{d,l}.
\end{equation*}

Несложно проверить (см. [4]), что при $ d \in \N, l \in \Z_+^d $
для $ \delta \in \R_+^d $ и $ x^0 \in \R^d $ таких, что $ x^0 +\delta I^d
\subset I^d, $ для $ f \in C(\overline I^d) $ и $ x \in \overline I^d $
выполняется равенство
\begin{equation*} \tag{2.1.9}
((\prod_{j=1}^d V_j(\mathcal P_{\delta_j, x_j^0}^{1, l_j}))f)(x) =
(\mathcal P_{\delta,x^0}^{d,l}f)(x).
\end{equation*}

В [4] установлена

Лемма 2.1.2

Пусть $ d \in \N, l \in \Z_+^d $ и $ \rho, \sigma \in \R_+^d. $
Тогда существует константа $ c_3(d,l,\rho,\sigma) >0 $ такая, что
для любых $ x^0, y^0 \in \R^d $ и $ \delta \in \R_+^d, $ для
которых $ (x^0 +\sigma \delta I^d) \subset (y^0 +\rho \delta I^d)
\subset I^d, $ для любого множества $ J \subset \Nu_{1,d}^1 $ для
$ f \in C(\overline I^d) $ соблюдается неравенство
\begin{equation*} \tag{2.1.10}
\|(\prod_{j \in J} V_j(E -\mathcal P_{\sigma_j
\delta_j,x_j^0}^{1,l_j}))f \|_{L_\infty(D)} \le c_3
\|f\|_{L_\infty(D)},
\end{equation*}
где $ D = y^0 +\rho \delta I^d.$
Определим операторы, с которыми мы будем иметь дело в следующем пункте, и
отметим их свойства, полезные для нас.

Для $ d \in \N, l,\kappa \in \Z_+^d, \nu \in \Nu_{0, 2^\kappa
-\e}^d $ определим непрерывный линейный оператор $  \mathtt R_{\kappa,
\nu}^{d,l}: C(\overline I^d) \mapsto \mathcal P^{d,l} \cap
C(\overline I^d) $ равенством
$$
\mathtt R_{\kappa, \nu}^{d,l} f =
\mathcal P_{\delta, x^0}^{d,l} (f \mid_{(x^0 +\delta \overline I^d)}),
f \in C(\overline I^d),
$$
при $ \delta = 2^{-\kappa}, x^0 = 2^{-\kappa} \nu. $

Как видно из (2.1.9), для $ d \in \N, l, \kappa \in \Z_+^d, \nu \in
\Nu_{0, 2^\kappa -\e}^d $ верно равенство
$$
 \mathtt R_{\kappa, \nu}^{d,l} =
\prod_{j=1}^d
V_j( \mathtt R_{\kappa_j, \nu_j}^{1,l_j}).
$$

При $ d \in \N $ для $ x \in \R^d $ положим
$$
x_+ = ((x_1)_+, \ldots, (x_d)_+),
$$
где $ t_+ = \frac{1} {2} (t +|t|), t \in \R. $

При $ d \in \N, l,m, \kappa \in \Z_+^d $
определим линейный оператор $ R_\kappa^{d,l,m}:  C(\overline I^d) \mapsto
 L_\infty (\R^d) $ равенством
$$
R_\kappa^{d,l,m} f = \sum_{ \nu \in \Nu_{-m, 2^\kappa -\e}^d}
(  \mathtt R_{\kappa, \nu_+}^{d,l} f)
g_{\kappa, \nu}^{d,m}.
$$

Как показано в [4], справедлива

Лемма 2.1.3

Пусть $ d \in \N, \alpha \in \R_+^d, 1 \le p < \infty $
удовлетворяют (1.3.7) и $ l = l(\alpha), m \in \Z_+^d. $ Тогда для любой
функции $ f \in \mathcal S_p^\alpha \mathcal H(I^d) $ справедливо соотношение
\begin{equation*} \tag{2.1.11}
\lim_{ \mn(\kappa) \to \infty}
\|f -R_\kappa^{d,l -\e,m}f\|_{L_p(I^d)} =0.
\end{equation*}

При $ d \in \N, l,m,\kappa \in \Z_+^d $ определим линейный оператор
$ \mathcal R_\kappa^{d,l,m}: C(\overline I^d) \mapsto L_\infty(\R^d),$ полагая
$$
\mathcal R_\kappa^{d,l,m} =
\sum_{\upsilon \in \Upsilon^d:  \s(\upsilon) \subset \s(\kappa)}
(-\e)^\upsilon R_{\kappa -\upsilon}^{d,l,m}.
$$

В [4] получено представление
\begin{equation*} \tag{2.1.12}
\mathcal R_\kappa^{d,l,m} =
\sum_{\nu \in \Nu_{-m, 2^\kappa -\e}^d}
g_{\kappa, \nu}^{d,m}
\mathcal U_{\kappa, \nu}^{d,l,m},
\end{equation*}
где $ \mathcal U_{\kappa, \nu}^{d,l,m}: C(\overline I^d) \mapsto \mathcal P^{d,l} $ --
линейный оператор, значение которого для $ f \in C(\overline I^d) $
определяется равенством
$$
\mathcal U_{\kappa, \nu}^{d,l,m} f =
\sum_{\upsilon \in \Upsilon^d: \s(\upsilon) \subset \s(\kappa)}
(-\e)^\upsilon
\sum_{\rho \in \Rho_{\kappa,\nu,\upsilon}^{d,m}}
(\prod_{j \in \s(\upsilon)}
a_{\nu_j -2\rho_j}^{m_j})
( \mathtt R_{\kappa -\upsilon, \rho_+}^{d,l} f),
$$
представляется в виде
\begin{multline*} \tag{2.1.13}
\mathcal U_{\kappa,\nu}^{d,l,m} =
\sum_{ \gamma \in \Upsilon^d: \s(\kappa) \subset \s(\gamma)}
(-\e)^\gamma
\sum_{\upsilon \in \Upsilon^d: \s(\upsilon) \subset \s(\kappa)}
(-\e)^\upsilon
\sum_{ \rho \in \Rho_{\kappa,\nu,\upsilon}^{d,m}}
(\prod_{j \in \s(\upsilon)}
a_{\nu_j -2\rho_j}^{m_j})\\
\times\biggl(\prod_{j \in \s(\gamma)}
V_j(E - \mathtt R_{\kappa_j -\upsilon_j, \rho_{+ j}}^{1,l_j})\biggr),
\end{multline*}
а множество
\begin{multline*}
\Rho_{\kappa,\nu,\upsilon}^{d,m} = \{ \rho \in \Nu_{-m, 2^{\kappa
-\upsilon} -\e}^d: \rho_j = \nu_j, j \in \Nu_{1,d}^1 \setminus
\s(\upsilon),  \\
(\nu_j -2 \rho_j) \in \Nu_{0, m_j +1}^1, j \in \s(\upsilon) \},
d \in \N, l,m,\kappa \in \Z_+^d, \\
\nu \in \Nu_{-m, 2^\kappa -\e}^d,
\upsilon \in \Upsilon^d: \s(\upsilon) \subset \s(\kappa).
\end{multline*}
\bigskip

2.2. В этом пункте будет установлена оценка сверху
наилучшей точности восстановления в $ L_q(I^d) $ производной
$ \D^\lambda f $ по значениям в $ n $ точках функций
$ f $ из классов $ \mathcal S_{p,\theta}^\alpha \mathcal  B(I^d). $

Сначала опишем постановку задачи.

Пусть $ T $ -- топологическое пространство, $ X $ --
банахово пространство над $ \R, U: D(U) \mapsto X $ --
линейный оператор с областью определения $ D(U) \subset C(T), $
принимающий значения в $ X. $
Пусть ещ\"е $ K \subset D(U) $ -- некоторый класс
функций.

Для $ n \in \N $ через $ \Phi_n(C(T)) $ обозначим совокупнсть
всех отображений $ \phi: C(T) \mapsto \R^n, $ для каждого
из которых существует набор точек $ \{ t^j \in T, j=1,\ldots,n\} $ такой, что
$$
\phi(f) = (f(t^1), \ldots, f(t^n)), f \in C(T).
$$

При $ n \in \N $ обозначим также через $ \mathcal A^n(X)
(\overline {\mathcal A}^n(X)) $ множество всех отображений (всех линейных отображений)
$ A: \R^n \mapsto X. $

Тогда при $ n \in \N $ положим
$$
\sigma_n(U,K,X) =
\inf_{ A \in \mathcal A^n(X), \phi \in \Phi_n(C(T))} \sup_{f \in K}
 \|U f -A \circ \phi (f)\|_X,
$$
а
$$
\overline \sigma_n(U,K,X) =
\inf_{ A \in \overline{\mathcal A}^n(X), \phi \in \Phi_n(C(T))} \sup_{f \in K}
 \|U f -A \circ \phi (f)\|_X.
$$

Лемму 2.2.1 можно рассматривать как модификацию леммы 3.4 из [4].

Лемма 2.2.1

Пусть $ d \in \N, \alpha \in \R_+^d, 1\le p < \infty, 1 \le q \le \infty, m \in
\Z_+^d, \lambda \in \Z_+^d(m) $ и выполняется (1.3.7). Тогда существуют
константы $ c_{1}(d,\alpha,p,q,m,\lambda) >0, c_{2}(d,m) >0,
c_{3}(d,\alpha) >0 $
и $ \epsilon \in \R_+^d $
такие, что для любой функции $ f \in S_p^\alpha H(I^d) $
при $ l = l(\alpha) $ и $ \kappa \in \Z_+^d \setminus \{0\} $
соблюдается неравенство
\begin{multline*} \tag{2.2.1}
\| \D^\lambda \mathcal R_\kappa^{d,l -\e,m} f
\|_{L_q(I^d)} \le c_{1}2^{(\kappa, \lambda +(p^{-1} -q^{-1})_+ \e)}
\biggl(\Omega^{l \chi_{\s(\kappa)}}(f, (c_{2} 2^{-\kappa})^{\s(\kappa)})_{L_p(I^d)}\\
+\sum_{J \subset \Nu_{1,d}^1: J \ne \emptyset}
\biggl(\int_{ (c_{3} I^d)^J}(\prod_{j \in J}u_j^{-p(\alpha_j -\epsilon_j) -1})
(\Omega^{l \chi_{J \cup \s(\kappa)}}( f, (u \chi_{J \setminus \s(\kappa)} \\
+2^{-\kappa} u \chi_{\s(\kappa) \cap J} +c_{2} 2^{-\kappa}
\chi_{\s(\kappa) \setminus J})^{J \cup \s(\kappa)})_{L_p(I^d)})^p du^J\biggr)^{1/p}\biggr).
\end{multline*}

Доказательство.

Сначала построим некоторые вспомогательные объекты, которые нам
понадобятся при доказательстве.

При $ d \in \N, m,\kappa \in \Z_+^d, \nu \in \Nu_{0, 2^\kappa
-\e}^d $ зададим точку $ \mathtt x_{\kappa, \nu}^{d,m} $ и вектор $
\mathtt \delta_{\kappa, \nu}^{d,m}, $
 полагая
\begin{multline*}
( \mathtt x_{\kappa, \nu}^{d,m})_j = 2^{-\kappa_j} \min((\nu_j
-2m_j -1)_+, (2^{\kappa_j} -2m_j -3)_+),
 (\mathtt \delta_{\kappa, \nu}^{d,m})_j \\
 =2^{-\kappa_j} \min(2^{\kappa_j}, 2m_j +3), j=1,\ldots,d,
\end{multline*} и определим клетку $ \mathtt D_{\kappa, \nu}^{d,m} $
равенством $ \mathtt D_{\kappa, \nu}^{d,m} =  \mathtt x_{\kappa,
\nu}^{d,m} + \mathtt \delta_{\kappa, \nu}^{d,m} I^d,$
а через $ \chi_{\mathtt D_{\kappa, \nu}^{d,m}} $ обозначим
характеристическую функцию множества $ \mathtt D_{\kappa, \nu}^{d,m}. $

Из определения следует, что при $ d \in \N, m,\kappa \in \Z_+^d,
\nu \in \Nu_{0, 2^\kappa -\e}^d $ множество
\begin{equation*} \tag{2.2.2}
 \mathtt D_{\kappa, \nu}^{d,m} \subset I^d,
\end{equation*}
а также
\begin{equation*} \tag{2.2.3}
Q_{\kappa,\nu}^d \subset
 \mathtt D_{\kappa, \nu}^{d,m}.
 \end{equation*}

С помощью (2.2.3) легко проверить, что при $ d \in \N, m \in
\Z_+^d $ существует константа $ c_{4}(d,m) >0 $ такая, что для $
\kappa \in \Z_+^d $ и $ x \in I^d $ верно неравенство
\begin{equation*} \tag{2.2.4}
\card \{ \nu \in \Nu_{0, 2^\kappa -\e}^d: x \in \mathtt D_{\kappa,
\nu}^{d,m} \} \le c_{4}(d,m).
\end{equation*}

Пользуясь тем, что для $ t \in \R, a \in \R_+ $ справедливы
соотношения $ (at)_+ = at_+ $ и $ t_+ \le (t+a)_+ \le t_+ +a, $ а
также используя (2.1.2), несложно показать, что при $ d \in \N,
m,\kappa \in \Z_+^d, $ для $ \nu^\prime \in \Nu_{0, 2^\kappa
-\e}^d, \nu \in \Nu_{-m, 2^\kappa -\e}^d: Q_{\kappa,
\nu^\prime}^d \cap \supp g_{\kappa, \nu}^{d,m} \ne \emptyset,
\upsilon \in \Upsilon^d: \s(\upsilon) \subset \s(\kappa), \rho \in
\Rho_{\kappa, \nu, \upsilon}^{d,m} $ имеет место включение
\begin{equation*} \tag{2.2.5}
Q_{\kappa -\upsilon, \rho_+}^d \subset \mathtt D_{\kappa,
\nu^\prime}^{d,m}.
\end{equation*}

Отметим ещ\"е, что
\begin{equation*} \tag{2.2.6}
\card \Rho_{\kappa,\nu,\upsilon}^{d,m} \le c_{5}(d,m),
\end{equation*}
$ d \in \N, m,\kappa \in \Z_+^d, \nu \in \Nu_{-m, 2^\kappa -\e}^d,
\upsilon \in \Upsilon^d: \s(\upsilon) \subset \s(\kappa). $

Определим ещ\"е при $ l, m, \kappa \in \Z_+, \nu \in \Nu_{0,
2^\kappa -1}^1 $ непрерывный линейный оператор $ \mathtt S_{\kappa,
\nu}^{ l, m}: L_1(I) \mapsto \mathcal P^{1, l}(I) \cap
L_\infty (I), $ полагая для $ f \in L_1(I) $ значение
$$
\mathtt S_{\kappa, \nu}^{ l, m}f = \mathtt P_{\delta, x^0}^{1, l}
(f\mid_{(x^0 +\delta I)})
$$
при $ \delta = 2^{-\kappa} \min(2^\kappa, 2 m +3), x^0 =
2^{-\kappa} \min((\nu -2 m -1)_+, (2^\kappa -2 m -3)_+). $

Теперь в условиях леммы 2.2.1 в силу (2.1.5) и (2.1.12) так же,
как (1.55) в [4], получаем
\begin{equation*} \tag{2.2.7}
\| \D^\lambda \mathcal R_\kappa^{d, l-\e,m}f \|_{L_q(I^d)} \le
\sum_{\mu \in \Z_+^d(\lambda)} C_\lambda^\mu
\| \sum_{ \nu \in \Nu_{-m, 2^\kappa -\e}^d}
\D^\mu
(\mathcal U_{\kappa, \nu}^{d, l-\e,m}f)
\D^{\lambda -\mu}
g_{\kappa,\nu}^{d,m}
\|_{L_q(I^d)}.
\end{equation*}

Оценивая правую часть (2.2.7), подобно (1.56) в [4], при $ \mu \in
\Z_+^d(\lambda) $ выводим
\begin{multline*} \tag{2.2.8}
\biggl\| \sum_{ \nu \in \Nu_{-m, 2^\kappa -\e}^d}
\D^\mu
(\mathcal U_{\kappa, \nu}^{d, l-\e,m}f)
\D^{\lambda -\mu}
g_{\kappa,\nu}^{d,m}
\biggr\|_{L_q(I^d)}^q \le \\
\sum_{\nu^\prime \in \Nu_{0, 2^\kappa -\e}^d}
\biggl( \sum_{ \nu \in \Nu_{-m, 2^\kappa -\e}^d:
\supp g_{\kappa, \nu}^{d,m} \cap Q_{\kappa, \nu^\prime}^d \ne \emptyset}
\|\D^\mu
(\mathcal U_{\kappa, \nu}^{d, l-\e,m}f)
\D^{\lambda -\mu}
g_{\kappa,\nu}^{d,m}
\|_{L_q( Q_{\kappa, \nu^\prime}^d)}\biggr)^q.
\end{multline*}

Применяя (2.1.4), (1.1.1) и (2.2.3), находим, что при
$ \nu^\prime \in \Nu_{0, 2^\kappa -\e}^d,
\nu \in \Nu_{-m, 2^\kappa -\e}^d:
\supp g_{\kappa, \nu}^{d,m} \cap Q_{\kappa, \nu^\prime}^d \ne \emptyset, $
справедливо неравенство
\begin{multline*} \tag{2.2.9}
\|\D^\mu (\mathcal U_{\kappa, \nu}^{d, l-\e,m}f) \D^{\lambda -\mu}
g_{\kappa,\nu}^{d,m} \|_{L_q( Q_{\kappa, \nu^\prime}^d)} \\
\le \|\D^{\lambda -\mu} g_{\kappa,\nu}^{d,m} \|_{L_\infty(\R^d)}
\|\D^\mu (\mathcal U_{\kappa, \nu}^{d, l-\e,m}f) \|_{L_q(
Q_{\kappa, \nu^\prime}^d)} \\
\le c_{6} 2^{(\kappa, \lambda -\mu)}
\|\D^\mu (\mathcal U_{\kappa, \nu}^{d, l-\e,m}f) \|_{L_q(
Q_{\kappa, \nu^\prime}^d)} \\
\le c_{6} 2^{(\kappa, \lambda -\mu)}
c_{7} (2^{-\kappa})^{-\mu +q^{-1}\e} \| \mathcal U_{\kappa,
\nu}^{d, l-\e,m}f \|_{L_\infty( Q_{\kappa, \nu^\prime}^d)} \\
= c_{8}2^{(\kappa, \lambda -q^{-1}\e)} \| \mathcal U_{\kappa, \nu}^{d,
l-\e,m}f \|_{L_\infty( Q_{\kappa, \nu^\prime}^d)} \le c_{8}
2^{(\kappa, \lambda -q^{-1}\e)} \| \mathcal U_{\kappa, \nu}^{d,
l-\e,m}f \|_{L_\infty( \mathtt D_{\kappa, \nu^\prime}^{d,m})}.
\end{multline*}

Из (2.1.13) для
$ \nu^\prime \in \Nu_{0, 2^\kappa -\e}^d,
\nu \in \Nu_{-m, 2^\kappa -\e}^d:
\supp g_{\kappa, \nu}^{d,m} \cap Q_{\kappa, \nu^\prime}^d \ne \emptyset,  $
имеем неравенство
\begin{multline*} \tag{2.2.10}
\| \mathcal U_{\kappa, \nu}^{d, l-\e,m}f
\|_{L_\infty( \mathtt D_{\kappa, \nu^\prime}^{d,m})}\\
\le \sum_{ \gamma \in \Upsilon^d: \s(\kappa) \subset \s(\gamma)}
\sum_{\upsilon \in \Upsilon^d: \s(\upsilon) \subset \s(\kappa)}
\sum_{ \rho \in \Rho_{\kappa,\nu,\upsilon}^{d,m}}
(\prod_{j \in \s(\upsilon)}a_{\nu_j -2\rho_j}^{m_j})\\
\times\|(\prod_{j \in \s(\gamma)}
V_j(E - \mathtt R_{\kappa_j -\upsilon_j, \rho_{+ j}}^{1,l_j -1}))f
\|_{L_\infty( \mathtt D_{\kappa, \nu^\prime}^{d,m})}.
\end{multline*}

Далее, для
$ \nu^\prime \in \Nu_{0, 2^\kappa -\e}^d,
\nu \in \Nu_{-m, 2^\kappa -\e}^d:
\supp g_{\kappa, \nu}^{d,m} \cap Q_{\kappa, \nu^\prime}^d \ne \emptyset,
\gamma \in \Upsilon^d: \s(\kappa) \subset \s(\gamma),
\upsilon \in \Upsilon^d: \s(\upsilon) \subset \s(\kappa),
\rho \in \Rho_{\kappa,\nu,\upsilon}^{d,m}, $
используя (1.3.3), а затем с уч\"етом (2.2.5), (2.2.2)
применяя (2.1.10), получаем
\begin{multline*} \tag{2.2.11}
\|(\prod_{j \in \s(\gamma)}
V_j(E - \mathtt R_{\kappa_j -\upsilon_j, \rho_{+ j}}^{1,l_j -1}))f
\|_{L_\infty( \mathtt D_{\kappa, \nu^\prime}^{d,m})} \\
=\|(\prod_{j \in \s(\gamma) \setminus \s(\kappa)}
V_j(E - \mathtt R_{\kappa_j -\upsilon_j, \rho_{+ j}}^{1,l_j -1}))
((\prod_{j \in \s(\kappa)}
V_j(E - \mathtt R_{\kappa_j -\upsilon_j, \rho_{+ j}}^{1,l_j -1}))f)
\|_{L_\infty( \mathtt D_{\kappa, \nu^\prime}^{d,m})} \\
\le c_{9} \|(\prod_{j \in \s(\kappa)}
V_j(E - \mathtt R_{\kappa_j -\upsilon_j, \rho_{+ j}}^{1,l_j -1}))f
\|_{L_\infty( \mathtt D_{\kappa, \nu^\prime}^{d,m})},
\end{multline*}
и пользуясь тем, что в виду (2.1.8) при $ j =1,\ldots,d $
соблюдено равенство
$$
E - \mathtt R_{\kappa_j -\upsilon_j, \rho_{+ j}}^{1,l_j -1} =
(E - \mathtt R_{\kappa_j -\upsilon_j, \rho_{+ j}}^{1,l_j -1})
(E -\mathtt S_{\kappa_j, \nu_j^\prime}^{ l_j -1,m_j}),
$$
благодаря п. 2 леммы 1.3.2 и (1.3.3), а также в силу (2.2.5), (2.2.2),
(2.1.10) выводим
\begin{multline*} \tag{2.2.12}
\|(\prod_{j \in \s(\kappa)} V_j(E - \mathtt R_{\kappa_j -\upsilon_j,
\rho_{+ j}}^{1,l_j -1}))f \|_{L_\infty( \mathtt D_{\kappa,
\nu^\prime}^{d,m})}\\
 = \|(\prod_{j \in \s(\kappa)} V_j((E - \mathtt
R_{\kappa_j -\upsilon_j, \rho_{+ j}}^{1,l_j -1}) (E -\mathtt
S_{\kappa_j, \nu_j^\prime}^{ l_j -1,m_j})))f \|_{L_\infty( \mathtt
D_{\kappa, \nu^\prime}^{d,m})} \\
= \|(\prod_{j \in \s(\kappa)} (V_j(E
- \mathtt R_{\kappa_j -\upsilon_j, \rho_{+ j}}^{1,l_j -1}) V_j(E -\mathtt
S_{\kappa_j, \nu_j^\prime}^{ l_j -1,m_j})))f \|_{L_\infty( \mathtt
D_{\kappa, \nu^\prime}^{d,m})} \\
= \|(\prod_{j \in \s(\kappa)} V_j(E
- \mathtt R_{\kappa_j -\upsilon_j, \rho_{+ j}}^{1,l_j -1})) ((\prod_{j
\in \s(\kappa)} V_j(E -\mathtt S_{\kappa_j, \nu_j^\prime}^{ l_j
-1,m_j}))f) \|_{L_\infty( \mathtt D_{\kappa, \nu^\prime}^{d,m})} \le
\\
c_{10} \|(\prod_{j \in \s(\kappa)} V_j(E - \mathtt S_{\kappa_j,
\nu_j^\prime}^{ l_j -1,m_j}))f \|_{L_\infty( \mathtt D_{\kappa,
\nu^\prime}^{d,m})}.
\end{multline*}

Соединяя (2.2.10), (2.2.11), (2.2.12), и учитывая (2.2.6), находим, что для
$ \nu^\prime \in \Nu_{0, 2^\kappa -\e}^d,
\nu \in \Nu_{-m, 2^\kappa -\e}^d:
\supp g_{\kappa, \nu}^{d,m} \cap Q_{\kappa, \nu^\prime}^d \ne \emptyset, $
 выполняется неравенство
\begin{multline*} \tag{2.2.13}
\| \mathcal U_{\kappa, \nu}^{d, l-\e,m}f
\|_{L_\infty( \mathtt D_{\kappa, \nu^\prime}^{d,m})} \\
\le \sum_{ \gamma \in \Upsilon^d: \s(\kappa) \subset \s(\gamma)}
\sum_{\upsilon \in \Upsilon^d: \s(\upsilon) \subset \s(\kappa)}
\sum_{ \rho \in \Rho_{\kappa,\nu,\upsilon}^{d,m}}
c_{11}\|(\prod_{j \in \s(\kappa)}
V_j(E -\mathtt S_{\kappa_j, \nu_j^\prime}^{ l_j -1,m_j}))f
\|_{L_\infty( \mathtt D_{\kappa, \nu^\prime}^{d,m})} \\
\le c_{12}\|(\prod_{j \in \s(\kappa)}
V_j(E -\mathtt S_{\kappa_j, \nu_j^\prime}^{ l_j -1,m_j}))f
\|_{L_\infty( \mathtt D_{\kappa, \nu^\prime}^{d,m})}.
\end{multline*}

Применяя сначала к правой части (2.2.13) неравенство (1.3.6)
при $ \lambda =0, q = \infty, $ а затем используя (1.3.4), приходим к
неравенству
\begin{multline*} \tag{2.2.14}
\| \mathcal U_{\kappa, \nu}^{d, l-\e,m}f
\|_{L_\infty( \mathtt D_{\kappa, \nu^\prime}^{d,m})} \le
c_{13}
( \mathtt \delta_{\kappa, \nu^\prime}^{d,m})^{-p^{-1} \e}
\biggl(\|(\prod_{j \in \s(\kappa)}
V_j(E -\mathtt S_{\kappa_j, \nu_j^\prime}^{ l_j -1,m_j}))f
\|_{L_p( \mathtt D_{\kappa, \nu^\prime}^{d,m})} \\+
\sum_{J \subset \Nu_{1,d}^1: J \ne \emptyset}
(\prod_{j \in J}
( \mathtt \delta_{\kappa, \nu^\prime}^{d,m})_j^{p^{-1}})
\int_{ (c_{14}  \mathtt \delta_{\kappa, \nu^\prime}^{d,m} I^d)^J}
(\prod_{j \in J}t_j^{-2p^{-1} -1})\\
\times\biggl(\int_{ (t B^d)^J}
\int_{ ( \mathtt D_{\kappa, \nu^\prime}^{d,m})_\xi^{l \chi_J}}
|\Delta_\xi^{l \chi_J}((\prod_{j \in \s(\kappa)}
V_j(E -\mathtt S_{\kappa_j, \nu_j^\prime}^{ l_j -1,m_j}))f)(x)|^p
dx d\xi^J\biggr)^{1/p} dt^J\biggr) \\
\le c_{15}2^{(\kappa, \e) p^{-1}}\biggl(c_{16}(\prod_{j \in \s(\kappa)}
( \mathtt \delta_{\kappa, \nu^\prime}^{d,m})_j^{-1/p})
\biggl(\int_{ ( \mathtt \delta_{\kappa, \nu^\prime}^{d,m} B^d)^{\s(\kappa)}}
\int_{ ( \mathtt D_{\kappa, \nu^\prime}^{d,m})_\xi^{l \chi_{\s(\kappa)}}}
 |\Delta_\xi^{l \chi_{\s(\kappa)}} f(x)|^p
dx d\xi^{\s(\kappa)}\biggr)^{1/p} \\
+c_{17}\sum_{J \subset \Nu_{1,d}^1: J \ne \emptyset}
(\prod_{j \in J}2^{-\kappa_j p^{-1}})\int_{ (c_{3} 2^{-\kappa} I^d)^J}
(\prod_{j \in J}t_j^{-2p^{-1} -1})
\biggl(\int_{ (t B^d)^J}\biggl(c_{18}(\prod_{j \in \s(\kappa) \setminus J}
( \mathtt \delta_{\kappa, \nu^\prime}^{d,m})_j^{-1/p})\\
\times\biggl(\int_{ ( \mathtt \delta_{\kappa, \nu^\prime}^{d,m} B^d)^{\s(\kappa) \setminus J}}
\int_{ ( \mathtt D_{\kappa, \nu^\prime}^{d,m})_\xi^{l \chi_{J \cup \s(\kappa)}}}
 |\Delta_\xi^{l \chi_{J \cup \s(\kappa)}} f(x)|^p
dx d\xi^{\s(\kappa) \setminus J}\biggr)^{1/p}\biggr)^p d\xi^J\biggr)^{1/p} dt^J\biggr)\\
\le c_{19}2^{(\kappa, \e) p^{-1}}\biggl((\prod_{j \in \s(\kappa)}
2^{\kappa_j /p})\biggl(\int_{ (c_{2} 2^{-\kappa} B^d)^{\s(\kappa)}}
\int_{ ( \mathtt D_{\kappa, \nu^\prime}^{d,m})_\xi^{l \chi_{\s(\kappa)}}}
 |\Delta_\xi^{l \chi_{\s(\kappa)}} f(x)|^p
dx d\xi^{\s(\kappa)}\biggr)^{1/p} \\
+\sum_{J \subset \Nu_{1,d}^1: J \ne \emptyset}
(\prod_{j \in J}2^{-\kappa_j p^{-1}})(\prod_{j \in \s(\kappa) \setminus J}
2^{\kappa_j /p})\int_{ (c_{3} 2^{-\kappa} I^d)^J}(\prod_{j \in J}
t_j^{-2p^{-1} -1})\\
\times\biggl(\int_{ (t B^d)^J}\int_{ ( c_{2} 2^{-\kappa} B^d)^{\s(\kappa) \setminus J}}
\int_{ ( \mathtt D_{\kappa, \nu^\prime}^{d,m})_\xi^{l \chi_{J \cup \s(\kappa)}}}
 |\Delta_\xi^{l \chi_{J \cup \s(\kappa)}} f(x)|^p
dx d\xi^{\s(\kappa) \setminus J} d\xi^J\biggr)^{1/p} dt^J\biggr),
\end{multline*}
$ \nu^\prime \in \Nu_{0, 2^\kappa -\e}^d,
\nu \in \Nu_{-m, 2^\kappa -\e}^d:
\supp g_{\kappa, \nu}^{d,m} \cap Q_{\kappa, \nu^\prime}^d \ne \emptyset.  $

Фиксируя $ \epsilon \in \R_+^d $ так, чтобы соблюдалось условие
$ \alpha -p^{-1} \e -\epsilon >0, $ и применяя неравенство Г\"ельдера, для
$ J \subset \Nu_{1,d}^1: J \ne \emptyset,
\nu^\prime \in \Nu_{0, 2^\kappa -\e}^d $ имеем
\begin{multline*} \tag{2.2.15}
\int_{ (c_{3} 2^{-\kappa} I^d)^J}(\prod_{j \in J}t_j^{-2p^{-1} -1})\\
\times\biggl(\int_{ (t B^d)^J}
\int_{ ( c_{2} 2^{-\kappa} B^d)^{\s(\kappa) \setminus J}}
\int_{ ( \mathtt D_{\kappa, \nu^\prime}^{d,m})_\xi^{l \chi_{J \cup \s(\kappa)}}}
 |\Delta_\xi^{l \chi_{J \cup \s(\kappa)}} f(x)|^p
dx d\xi^{\s(\kappa) \setminus J} d\xi^J\biggr)^{1/p} dt^J \\
=\int_{ (c_{3} 2^{-\kappa} I^d)^J}(\prod_{j \in J}
t_j^{\alpha_j -p^{-1} -\epsilon_j -1/p^\prime})
(\prod_{j \in J}t_j^{-(\alpha_j -\epsilon_j +p^{-1}) -1/p})\\
\times\biggl(\int_{ (t B^d)^J}
\int_{ ( c_{2} 2^{-\kappa} B^d)^{\s(\kappa) \setminus J}}
\int_{ ( \mathtt D_{\kappa, \nu^\prime}^{d,m})_\xi^{l \chi_{J \cup \s(\kappa)}}}
 |\Delta_\xi^{l \chi_{J \cup \s(\kappa)}} f(x)|^p
dx d\xi^{\s(\kappa) \setminus J} d\xi^J\biggr)^{1/p} dt^J\\
 \le\biggl(\int_{ (c_{3} 2^{-\kappa} I^d)^J}
(\prod_{j \in J}
t_j^{p^\prime (\alpha_j -p^{-1} -\epsilon_j) -1})
dt^J\biggr)^{1/p^\prime}\biggl(\int_{ (c_{3} 2^{-\kappa} I^d)^J}
(\prod_{j \in J}t_j^{-p(\alpha_j -\epsilon_j +p^{-1}) -1})\\
\times\int_{ (t B^d)^J}\int_{ ( c_{2} 2^{-\kappa} B^d)^{\s(\kappa) \setminus J}}
\int_{ ( \mathtt D_{\kappa, \nu^\prime}^{d,m})_\xi^{l \chi_{J \cup \s(\kappa)}}}
 |\Delta_\xi^{l \chi_{J \cup \s(\kappa)}} f(x)|^p
dx d\xi^{\s(\kappa) \setminus J} d\xi^J dt^J\biggr)^{1/p} \\
=c_{20}(\prod_{j \in J}2^{-\kappa_j (\alpha_j -p^{-1} -\epsilon_j)})
\biggl(\int_{ (c_{3} 2^{-\kappa} I^d)^J}
(\prod_{j \in J}
t_j^{-p(\alpha_j -\epsilon_j +p^{-1}) -1})\\
\times\int_{ (t B^d)^J}\int_{ ( c_{2} 2^{-\kappa} B^d)^{\s(\kappa) \setminus J}}
\int_{ ( \mathtt D_{\kappa, \nu^\prime}^{d,m})_\xi^{l \chi_{J \cup \s(\kappa)}}}
 |\Delta_\xi^{l \chi_{J \cup \s(\kappa)}} f(x)|^p
dx d\xi^{\s(\kappa) \setminus J} d\xi^J dt^J\biggr)^{1/p}.
\end{multline*}

Соединяя (2.2.9), (2.2.14), (2.2.15), заключаем, что при
$ \nu^\prime \in \Nu_{0, 2^\kappa -\e}^d,
\nu \in \Nu_{-m, 2^\kappa -\e}^d:
\supp g_{\kappa, \nu}^{d,m} \cap Q_{\kappa, \nu^\prime}^d \ne \emptyset, $
справедливо неравенство
\begin{multline*} \tag{2.2.16}
\|\D^\mu
(\mathcal U_{\kappa, \nu}^{d, l-\e,m}f)
\D^{\lambda -\mu}
g_{\kappa,\nu}^{d,m}
\|_{L_q( Q_{\kappa, \nu^\prime}^d)}
\le c_{21}2^{(\kappa, \lambda +p^{-1} \e -q^{-1} \e)}\\
\times\biggl((\prod_{j \in \s(\kappa)}
2^{\kappa_j /p})
\biggl(\int_{ (c_{2} 2^{-\kappa} B^d)^{\s(\kappa)}}
\int_{ ( \mathtt D_{\kappa, \nu^\prime}^{d,m})_\xi^{l \chi_{\s(\kappa)}}}
 |\Delta_\xi^{l \chi_{\s(\kappa)}} f(x)|^p
dx d\xi^{\s(\kappa)}\biggr)^{1/p} \\
+\sum_{J \subset \Nu_{1,d}^1: J \ne \emptyset}
(\prod_{j \in J}
2^{-\kappa_j (\alpha_j -\epsilon_j)})
(\prod_{j \in \s(\kappa) \setminus J}
2^{\kappa_j /p})\biggl(\int_{ (c_{3} 2^{-\kappa} I^d)^J}
(\prod_{j \in J}
t_j^{-p(\alpha_j -\epsilon_j +p^{-1}) -1})\\
\times
\int_{ (t B^d)^J}
\int_{ ( c_{2} 2^{-\kappa} B^d)^{\s(\kappa) \setminus J}}
\int_{ ( \mathtt D_{\kappa, \nu^\prime}^{d,m})_\xi^{l \chi_{J \cup \s(\kappa)}}}
 |\Delta_\xi^{l \chi_{J \cup \s(\kappa)}} f(x)|^p
dx d\xi^{\s(\kappa) \setminus J} d\xi^J dt^J\biggr)^{1/p}\biggr).
\end{multline*}

Подставляя (2.2.16) в (2.2.8) и учитывая (2.1.3), а затем
применяя неравенства Г\"ельдера с показателями $ q $ и $ p/q \le 1, $
и используя (2.2.4), получаем, что при $ p \le q $
и $ \mu \in \Z_+^d(\lambda) $ выполняется неравенство
\begin{multline*}
\| \sum_{ \nu \in \Nu_{-m, 2^\kappa -\e}^d}
\D^\mu
(\mathcal U_{\kappa, \nu}^{d, l-\e,m}f)
\D^{\lambda -\mu}
g_{\kappa,\nu}^{d,m}\|_{L_q(I^d)}^q
 \le\sum_{\nu^\prime \in \Nu_{0, 2^\kappa -\e}^d}
\biggl( c_{22} 2^{(\kappa, \lambda +p^{-1} \e -q^{-1} \e)}\\
\times\biggl((\prod_{j \in \s(\kappa)}
2^{\kappa_j /p})
\biggl(\int_{ (c_{2} 2^{-\kappa} B^d)^{\s(\kappa)}}
\int_{ ( \mathtt D_{\kappa, \nu^\prime}^{d,m})_\xi^{l \chi_{\s(\kappa)}}}
 |\Delta_\xi^{l \chi_{\s(\kappa)}} f(x)|^p
dx d\xi^{\s(\kappa)}\biggr)^{1/p} \\
+
\sum_{J \subset \Nu_{1,d}^1: J \ne \emptyset}
(\prod_{j \in J}
2^{-\kappa_j (\alpha_j -\epsilon_j)})
(\prod_{j \in \s(\kappa) \setminus J}
2^{\kappa_j /p})
\biggl(\int_{ (c_{3} 2^{-\kappa} I^d)^J}
(\prod_{j \in J}
t_j^{-p(\alpha_j -\epsilon_j +p^{-1}) -1})\\
\times\int_{ (t B^d)^J}
\int_{ ( c_{2} 2^{-\kappa} B^d)^{\s(\kappa) \setminus J}}
\int_{ ( \mathtt D_{\kappa, \nu^\prime}^{d,m})_\xi^{l \chi_{J \cup \s(\kappa)}}}
 |\Delta_\xi^{l \chi_{J \cup \s(\kappa)}} f(x)|^p
dx d\xi^{\s(\kappa) \setminus J} d\xi^J dt^J\biggr)^{1/p}\biggr)\biggr)^q \\
\le\sum_{\nu^\prime \in \Nu_{0, 2^\kappa -\e}^d}
(c_{23}
2^{(\kappa, \lambda +p^{-1} \e -q^{-1} \e)})^q
\biggl((\prod_{j \in \s(\kappa)}
2^{\kappa_j q /p})\\
\times\biggl(\int_{ (c_{2} 2^{-\kappa} B^d)^{\s(\kappa)}}
\int_{ ( \mathtt D_{\kappa, \nu^\prime}^{d,m})_\xi^{l \chi_{\s(\kappa)}}}
 |\Delta_\xi^{l \chi_{\s(\kappa)}} f(x)|^p
dx d\xi^{\s(\kappa)}\biggr)^{q/p} \\
+\sum_{J \subset \Nu_{1,d}^1: J \ne \emptyset}
(\prod_{j \in J}
2^{-\kappa_j q(\alpha_j -\epsilon_j)})
(\prod_{j \in \s(\kappa) \setminus J}
2^{\kappa_j q /p})
\biggl(\int_{ (c_{3} 2^{-\kappa} I^d)^J}
(\prod_{j \in J}
t_j^{-p(\alpha_j -\epsilon_j +p^{-1}) -1})\\
\times\int_{ (t B^d)^J}
\int_{ ( c_{2} 2^{-\kappa} B^d)^{\s(\kappa) \setminus J}}
\int_{ ( \mathtt D_{\kappa, \nu^\prime}^{d,m})_\xi^{l \chi_{J \cup \s(\kappa)}}}
 |\Delta_\xi^{l \chi_{J \cup \s(\kappa)}} f(x)|^p
dx d\xi^{\s(\kappa) \setminus J} d\xi^J dt^J\biggr)^{q/p}\biggr) \\
=(c_{23}
2^{(\kappa, \lambda +p^{-1} \e -q^{-1} \e)})^q
\biggl((\prod_{j \in \s(\kappa)}
2^{\kappa_j q /p})\\
\times\sum_{\nu^\prime \in \Nu_{0, 2^\kappa -\e}^d}
\biggl(\int_{ (c_{2} 2^{-\kappa} B^d)^{\s(\kappa)}}
\int_{ ( \mathtt D_{\kappa, \nu^\prime}^{d,m})_\xi^{l \chi_{\s(\kappa)}}}
 |\Delta_\xi^{l \chi_{\s(\kappa)}} f(x)|^p
dx d\xi^{\s(\kappa)}\biggr)^{q/p}\\
 +\sum_{J \subset \Nu_{1,d}^1: J \ne \emptyset}
(\prod_{j \in J}
2^{-\kappa_j q(\alpha_j -\epsilon_j)})
(\prod_{j \in \s(\kappa) \setminus J}
2^{\kappa_j q /p})\\
\times\sum_{\nu^\prime \in \Nu_{0, 2^\kappa -\e}^d}
\biggl(\int_{ (c_{3} 2^{-\kappa} I^d)^J}
(\prod_{j \in J}
t_j^{-p(\alpha_j -\epsilon_j +p^{-1}) -1})\\
\times\int_{ (t B^d)^J}
\int_{ ( c_{2} 2^{-\kappa} B^d)^{\s(\kappa) \setminus J}}
\int_{ ( \mathtt D_{\kappa, \nu^\prime}^{d,m})_\xi^{l \chi_{J \cup \s(\kappa)}}}
 |\Delta_\xi^{l \chi_{J \cup \s(\kappa)}} f(x)|^p
dx d\xi^{\s(\kappa) \setminus J} d\xi^J dt^J\biggr)^{q/p}\biggr)
\displaybreak\\
 \le (c_{23}
2^{(\kappa, \lambda +p^{-1} \e -q^{-1} \e)})^q
\biggl((\prod_{j \in \s(\kappa)}
2^{\kappa_j q /p})\\
\times\biggl(\sum_{\nu^\prime \in \Nu_{0, 2^\kappa -\e}^d}
\int_{ (c_{2} 2^{-\kappa} B^d)^{\s(\kappa)}}
\int_{ ( \mathtt D_{\kappa, \nu^\prime}^{d,m})_\xi^{l \chi_{\s(\kappa)}}}
 |\Delta_\xi^{l \chi_{\s(\kappa)}} f(x)|^p
dx d\xi^{\s(\kappa)}\biggr)^{q/p}\\
 +\sum_{J \subset \Nu_{1,d}^1: J \ne \emptyset}
(\prod_{j \in J}
2^{-\kappa_j q(\alpha_j -\epsilon_j)})
(\prod_{j \in \s(\kappa) \setminus J}
2^{\kappa_j q /p})\\
\times\biggl(\sum_{\nu^\prime \in \Nu_{0, 2^\kappa -\e}^d}
\int_{ (c_{3} 2^{-\kappa} I^d)^J}
(\prod_{j \in J}
t_j^{-p(\alpha_j -\epsilon_j +p^{-1}) -1})\\
\times\int_{ (t B^d)^J}
\int_{ ( c_{2} 2^{-\kappa} B^d)^{\s(\kappa) \setminus J}}
\int_{ ( \mathtt D_{\kappa, \nu^\prime}^{d,m})_\xi^{l \chi_{J \cup \s(\kappa)}}}
 |\Delta_\xi^{l \chi_{J \cup \s(\kappa)}} f(x)|^p
dx d\xi^{\s(\kappa) \setminus J} d\xi^J dt^J\biggr)^{q/p}\biggr)\\
 \le (c_{23}
2^{(\kappa, \lambda +p^{-1} \e -q^{-1} \e)})^q
\biggl((\prod_{j \in \s(\kappa)}
2^{\kappa_j q /p})\\
\times\biggl(\sum_{\nu^\prime \in \Nu_{0, 2^\kappa -\e}^d}
\int_{ (c_{2} 2^{-\kappa} B^d)^{\s(\kappa)}}
\int_{ ( I^d)_\xi^{l \chi_{\s(\kappa)}}}
\chi_{\mathtt D_{\kappa, \nu^\prime}^{d,m}}(x)
 |\Delta_\xi^{l \chi_{\s(\kappa)}} f(x)|^p
dx d\xi^{\s(\kappa)}\biggr)^{q/p} \\
+\sum_{J \subset \Nu_{1,d}^1: J \ne \emptyset}
(\prod_{j \in J}
2^{-\kappa_j q(\alpha_j -\epsilon_j)})
(\prod_{j \in \s(\kappa) \setminus J}
2^{\kappa_j q /p})
\biggl(\sum_{\nu^\prime \in \Nu_{0, 2^\kappa -\e}^d}
\int_{ (c_{3} 2^{-\kappa} I^d)^J}
(\prod_{j \in J}
t_j^{-p(\alpha_j -\epsilon_j +p^{-1}) -1})\\
\times\int_{ (t B^d)^J}
\int_{ ( c_{2} 2^{-\kappa} B^d)^{\s(\kappa) \setminus J}}
\int_{ ( I^d)_\xi^{l \chi_{J \cup \s(\kappa)}}}
\chi_{\mathtt D_{\kappa, \nu^\prime}^{d,m}}(x)
 |\Delta_\xi^{l \chi_{J \cup \s(\kappa)}} f(x)|^p
dx d\xi^{\s(\kappa) \setminus J} d\xi^J dt^J\biggr)^{q/p}\biggr)\\
 =(c_{23}
2^{(\kappa, \lambda +p^{-1} \e -q^{-1} \e)})^q
\biggl((\prod_{j \in \s(\kappa)}
2^{\kappa_j q /p})\\
\times\biggl(\int_{ (c_{2} 2^{-\kappa} B^d)^{\s(\kappa)}}
\int_{ ( I^d)_\xi^{l \chi_{\s(\kappa)}}}
(\sum_{\nu^\prime \in \Nu_{0, 2^\kappa -\e}^d}
\chi_{\mathtt D_{\kappa, \nu^\prime}^{d,m}}(x))
 |\Delta_\xi^{l \chi_{\s(\kappa)}} f(x)|^p
dx d\xi^{\s(\kappa)}\biggr)^{q/p}\\
 +\sum_{J \subset \Nu_{1,d}^1: J \ne \emptyset}
(\prod_{j \in J}
2^{-\kappa_j q(\alpha_j -\epsilon_j)})
(\prod_{j \in \s(\kappa) \setminus J}
2^{\kappa_j q /p})
\biggl(\int_{ (c_{3} 2^{-\kappa} I^d)^J}
(\prod_{j \in J}
t_j^{-p(\alpha_j -\epsilon_j +p^{-1}) -1})\\
\times\int_{ (t B^d)^J}
\int_{ ( c_{2} 2^{-\kappa} B^d)^{\s(\kappa) \setminus J}}
\int_{ ( I^d)_\xi^{l \chi_{J \cup \s(\kappa)}}}
(\sum_{\nu^\prime \in \Nu_{0, 2^\kappa -\e}^d}
\chi_{\mathtt D_{\kappa, \nu^\prime}^{d,m}}(x))
 |\Delta_\xi^{l \chi_{J \cup \s(\kappa)}} f(x)|^p
dx d\xi^{\s(\kappa) \setminus J} d\xi^J dt^J\biggr)^{q/p}\biggr)\\
 \le (c_{23}
2^{(\kappa, \lambda +p^{-1} \e -q^{-1} \e)})^q
\biggl((\prod_{j \in \s(\kappa)}
2^{\kappa_j q /p})
\biggl(\int_{ (c_{2} 2^{-\kappa} B^d)^{\s(\kappa)}}
\int_{ ( I^d)_\xi^{l \chi_{\s(\kappa)}}}
c_{4}
 |\Delta_\xi^{l \chi_{\s(\kappa)}} f(x)|^p
dx d\xi^{\s(\kappa)}\biggr)^{q/p} \\
+\sum_{J \subset \Nu_{1,d}^1: J \ne \emptyset}
(\prod_{j \in J}
2^{-\kappa_j q(\alpha_j -\epsilon_j)})
(\prod_{j \in \s(\kappa) \setminus J}
2^{\kappa_j q /p})
\biggl(\int_{ (c_{3} 2^{-\kappa} I^d)^J}
(\prod_{j \in J}
t_j^{-p(\alpha_j -\epsilon_j +p^{-1}) -1})\\
\int_{ (t B^d)^J}
\int_{ ( c_{2} 2^{-\kappa} B^d)^{\s(\kappa) \setminus J}}
\int_{ ( I^d)_\xi^{l \chi_{J \cup \s(\kappa)}}}
c_{4}
 |\Delta_\xi^{l \chi_{J \cup \s(\kappa)}} f(x)|^p
dx d\xi^{\s(\kappa) \setminus J} d\xi^J dt^J\biggr)^{q/p}\biggr),
\end{multline*}
или
\begin{multline*} \tag{2.2.17}
\| \sum_{ \nu \in \Nu_{-m, 2^\kappa -\e}^d}
\D^\mu
(\mathcal U_{\kappa, \nu}^{d, l-\e,m}f)
\D^{\lambda -\mu}
g_{\kappa,\nu}^{d,m}\|_{L_q(I^d)} \\
\le c_{24}
2^{(\kappa, \lambda +p^{-1} \e -q^{-1} \e)}
\biggl((\prod_{j \in \s(\kappa)}
2^{\kappa_j q /p})
\biggl(\int_{ (c_{2} 2^{-\kappa} B^d)^{\s(\kappa)}}
\int_{ ( I^d)_\xi^{l \chi_{\s(\kappa)}}}
 |\Delta_\xi^{l \chi_{\s(\kappa)}} f(x)|^p
dx d\xi^{\s(\kappa)}\biggr)^{q/p} \\
+\sum_{J \subset \Nu_{1,d}^1: J \ne \emptyset}
(\prod_{j \in J}
2^{-\kappa_j q(\alpha_j -\epsilon_j)})
(\prod_{j \in \s(\kappa) \setminus J}
2^{\kappa_j q /p})
\biggl(\int_{ (c_{3} 2^{-\kappa} I^d)^J}
(\prod_{j \in J}
t_j^{-p(\alpha_j -\epsilon_j +p^{-1}) -1})\\
\times\int_{ (t B^d)^J}
\int_{ ( c_{2} 2^{-\kappa} B^d)^{\s(\kappa) \setminus J}}
\int_{ ( I^d)_\xi^{l \chi_{J \cup \s(\kappa)}}}
 |\Delta_\xi^{l \chi_{J \cup \s(\kappa)}} f(x)|^p
dx d\xi^{\s(\kappa) \setminus J} d\xi^J dt^J\biggr)^{q/p} \biggr)^{1/q}\\
 \le c_{24}
2^{(\kappa, \lambda +p^{-1} \e -q^{-1} \e)}
\biggl((\prod_{j \in \s(\kappa)}
2^{\kappa_j /p})
\biggl(\int_{ (c_{2} 2^{-\kappa} B^d)^{\s(\kappa)}}
\int_{ ( I^d)_\xi^{l \chi_{\s(\kappa)}}}
 |\Delta_\xi^{l \chi_{\s(\kappa)}} f(x)|^p
dx d\xi^{\s(\kappa)}\biggr)^{1/p} \\
+\sum_{J \subset \Nu_{1,d}^1: J \ne \emptyset}
(\prod_{j \in J}
2^{-\kappa_j (\alpha_j -\epsilon_j)})
(\prod_{j \in \s(\kappa) \setminus J}
2^{\kappa_j /p})
\biggl(\int_{ (c_{3} 2^{-\kappa} I^d)^J}
(\prod_{j \in J}
t_j^{-p(\alpha_j -\epsilon_j +p^{-1}) -1})\\
\times\int_{ (t B^d)^J}
\int_{ ( c_{2} 2^{-\kappa} B^d)^{\s(\kappa) \setminus J}}
\int_{ ( I^d)_\xi^{l \chi_{J \cup \s(\kappa)}}}
 |\Delta_\xi^{l \chi_{J \cup \s(\kappa)}} f(x)|^p
dx d\xi^{\s(\kappa) \setminus J} d\xi^J dt^J\biggr)^{1/p}\biggr)\\
\le c_{24}
2^{(\kappa, \lambda +p^{-1} \e -q^{-1} \e)}
\biggl((\prod_{j \in \s(\kappa)}
2^{\kappa_j /p})
((\prod_{j \in \s(\kappa)} (2 c_{2} 2^{-\kappa_j}))\\
\times\supvrai_{ \xi \in \R^d: \xi^{\s(\kappa)} \in
(c_{2} 2^{-\kappa} B^d)^{\s(\kappa)}}
\|\Delta_\xi^{l \chi_{\s(\kappa)}} f \|_{L_p
(( I^d)_\xi^{l \chi_{\s(\kappa)}})}^p )^{1/p} \\
+\sum_{J \subset \Nu_{1,d}^1: J \ne \emptyset}
(\prod_{j \in J}
2^{-\kappa_j (\alpha_j -\epsilon_j)})
(\prod_{j \in \s(\kappa) \setminus J}
2^{\kappa_j /p})
\biggl(\int_{ (c_{3} 2^{-\kappa} I^d)^J}
(\prod_{j \in J}
t_j^{-p(\alpha_j -\epsilon_j +p^{-1}) -1})
\biggl((\prod_{j \in J} (2 t_j))\\
\times
(\prod_{j \in (\s(\kappa) \setminus J)} ( 2 c_{2} 2^{-\kappa_j}))
\supvrai_{ \xi \in \R^d: \xi^J \in (t B^d)^J,
\xi^{\s(\kappa) \setminus J} \in
( c_{2} 2^{-\kappa} B^d)^{\s(\kappa) \setminus J}}
\| \Delta_\xi^{l \chi_{J \cup \s(\kappa)}} f \|_{L_p
(( I^d)_\xi^{l \chi_{J \cup \s(\kappa)}})}^p \biggr) dt^J\biggr)^{1/p}\biggr)\\
 \le c_{25}
2^{(\kappa, \lambda +p^{-1} \e -q^{-1} \e)}
\biggl(\Omega^{l \chi_{\s(\kappa)}}(f, (c_{2} 2^{-\kappa})^{\s(\kappa)})_{L_p(I^d)}\\
+\sum_{J \subset \Nu_{1,d}^1: J \ne \emptyset}
(\prod_{j \in J}
2^{-\kappa_j (\alpha_j -\epsilon_j)})
\biggl(\int_{ (c_{3} 2^{-\kappa} I^d)^J}
(\prod_{j \in J}
t_j^{-p(\alpha_j -\epsilon_j) -1})\\
(\Omega^{l \chi_{J \cup \s(\kappa)}}( f,
(t \chi_J +c_{2} 2^{-\kappa}
\chi_{\s(\kappa) \setminus J})^{J \cup \s(\kappa)})_{L_p(I^d)})^p dt^J\biggr)^{1/p}\biggr)
\displaybreak\\
= c_{25}
2^{(\kappa, \lambda +p^{-1} \e -q^{-1} \e)}
\biggl(\Omega^{l \chi_{\s(\kappa)}}(f, (c_{2} 2^{-\kappa})^{\s(\kappa)})_{L_p(I^d)}\\
+\sum_{J \subset \Nu_{1,d}^1: J \ne \emptyset}
(\prod_{j \in J}
2^{-\kappa_j (\alpha_j -\epsilon_j)})
\biggl(\int_{ (c_{3} I^d)^J}
(\prod_{j \in J}
(2^{-\kappa_j} u_j)^{-p(\alpha_j -\epsilon_j) -1})\\
\times(\Omega^{l \chi_{J \cup \s(\kappa)}}( f,
(2^{-\kappa} u \chi_J +c_{2} 2^{-\kappa}
\chi_{\s(\kappa) \setminus J})^{J \cup \s(\kappa)})_{L_p(I^d)})^p
(\prod_{j \in J} 2^{-\kappa_j}) du^J\biggr)^{1/p}\biggr)\\
=c_{25}
2^{(\kappa, \lambda +p^{-1} \e -q^{-1} \e)}
\biggl(\Omega^{l \chi_{\s(\kappa)}}(f, (c_{2} 2^{-\kappa})^{\s(\kappa)})_{L_p(I^d)}\\
+\sum_{J \subset \Nu_{1,d}^1: J \ne \emptyset}
\biggl(\int_{ (c_{3} I^d)^J}
(\prod_{j \in J}
u_j^{-p(\alpha_j -\epsilon_j) -1})
(\Omega^{l \chi_{J \cup \s(\kappa)}}( f, (u \chi_{J \setminus \s(\kappa)} \\
+2^{-\kappa} u \chi_{\s(\kappa) \cap J} +c_{2} 2^{-\kappa}
\chi_{\s(\kappa) \setminus J})^{J \cup \s(\kappa)})_{L_p(I^d)})^p du^J\biggr)^{1/p}\biggr).
\end{multline*}

Подставляя (2.2.17) в (2.2.7), находим, что при $ p \le q $ и $ \kappa
\in \Z_+^d \setminus \{0\} $ для  $ f \in S_p^\alpha H(I^d) $
соблюдается неравенство
\begin{multline*}
\| \D^\lambda \mathcal R_\kappa^{d, l-\e,m}f \|_{L_q(I^d)} \le
c_{26}
2^{(\kappa, \lambda +p^{-1} \e -q^{-1} \e)}
\biggl(\Omega^{l \chi_{\s(\kappa)}}(f, (c_{2} 2^{-\kappa})^{\s(\kappa)})_{L_p(I^d)}\\
+\sum_{J \subset \Nu_{1,d}^1: J \ne \emptyset}
\biggl(\int_{ (c_{3} I^d)^J}
(\prod_{j \in J}
u_j^{-p(\alpha_j -\epsilon_j) -1})
(\Omega^{l \chi_{J \cup \s(\kappa)}}( f, (u \chi_{J \setminus \s(\kappa)} \\
+2^{-\kappa} u \chi_{\s(\kappa) \cap J} +c_{2} 2^{-\kappa}
\chi_{\s(\kappa) \setminus J})^{J \cup \s(\kappa)})_{L_p(I^d)})^p du^J\biggr)^{1/p}\biggr),
\end{multline*}
которое совпадает с (2.2.1) при $ p \le q. $
Неравенство (2.2.1) при $ q < p $ вытекает из неравенства
(2.2.1) при $ q = p $ и того факта, что для
$ h \in L_p(I^d) $ при $ q < p $ справелдива оценка
$ \|h\|_{L_q(I^d)} \le  \|h\|_{L_p(I^d)}. \square $

Предложение 2.2.2

Пусть $ d \in \N, \alpha \in \R_+^d, l =  l(\alpha),
1 \le p < \infty, 1 \le q \le \infty, m \in \Z_+^d,
\lambda \in \Z_+^d(m) $ и выполняются условия (1.3.5), (1.3.7).
Тогда для $ f \in \mathcal S_p^\alpha \mathcal H(I^d) $
в $ L_q(I^d) $ имеет место равенство
\begin{equation*} \tag{2.2.18}
\D^\lambda f = \sum_{\kappa \in \Z_+^d}
\D^\lambda \mathcal R_\kappa^{d, l -\e,m}f.
\end{equation*}

Предложение 2.2.2 имеется в [4], однако там сходимость ряда, стоящего
в правой части (2.2.18), понимается несколько в ином смысле, чем в настоящей
работе. Поэтому ниже приводится доказательство предложения 2.2.2.

Доказательство.

В условиях предложения из (2.1.11) на основании леммы 1.2.1 заключаем, что для
$ f \in \mathcal S_p^\alpha \mathcal H(I^d) $
в $ L_p(I^d) $ справедливо равенство
\begin{equation*} \tag{2.2.19}
f = \sum_{ \kappa \in \Z_+^d}
\mathcal R_\kappa^{d, l -\e,m}f.
\end{equation*}

Из (2.2.1) и (1.3.5) (см. также (3.19) в [4]) следует, что для $ f
\in \mathcal S_p^\alpha \mathcal H(I^d) $ ряд
$$
\sum_{\kappa \in \Z_+^d}
\| \D^\lambda \mathcal R_\kappa^{d, l -\e,m}f
\|_{L_q(I^d)}
$$
сходится, и, значит, ряд
$$
\sum_{\kappa \in \Z_+^d}
 \D^\lambda \mathcal R_\kappa^{d, l -\e,m}f
$$
сходится в $ L_q(I^d). $

Принимая во внимание это обстоятельство и равенство (2.2.19),
для любой бесконечно дифференцируемой финитной функции $ \phi,$
носитель которой $ \supp \phi \subset I^d, $ имеем
\begin{multline*}
\langle \D^\lambda f, \phi \rangle =
(-1)^{(\lambda,\e)}
\langle f, \D^\lambda \phi \rangle =
(-1)^{(\lambda,\e)}
\int_{I^d} f \D^\lambda \phi dx \\
 =\int_{I^d}
(\sum_{ \kappa \in \Z_+^d} \mathcal R_\kappa^{d, l -\e,m}f)
(-1)^{(\lambda,\e)}
\D^\lambda \phi dx =
\sum_{ \kappa \in \Z_+^d}
(-1)^{(\lambda,\e)}
\int_{I^d}(\mathcal R_\kappa^{d, l -\e,m}f)
\D^\lambda \phi dx \\
=\sum_{ \kappa \in \Z_+^d}
\int_{I^d}
\D^\lambda(\mathcal R_\kappa^{d, l -\e,m}f)
\phi dx =
\int_{I^d}
(\sum_{ \kappa \in \Z_+^d}
\D^\lambda(\mathcal R_\kappa^{d, l -\e,m}f))
\phi dx.
\end{multline*}
А это означает, что в $ L_q(I^d) $ имеет место равенство (2.2.18). $ \square $

Теорема 2.2.3

Пусть $ d \in \N, \alpha \in \R_+^d, 1\le p < \infty, 1 \le q \le \infty,
\lambda \in \Z_+^d $ удовлетворяют условиям (1.3.5), (1.3.7).
Пусть ещ\"е $ U = \D^\lambda, D(U) = \{f \in C(I^d): \D^\lambda f \in L_q(I^d) \},
X = L_q(I^d), K = \mathcal S_{p,\theta}^\alpha \mathcal B(I^d), 1 \le \theta
\le \infty. $
Тогда существует константа $ c_{27}(d,\alpha,p,q,\theta,\lambda) >0 $ и
$ n_0 \in \N $ такие,что при $ n \ge n_0 $ выполняется неравенство
\begin{equation*} \tag{2.2.20}
\overline \sigma_n(U,K,X) \le c_{27} n^{-\mn}
(\log n)^{(\mn +1 -1/\max(p,\theta))(\cmn -1)},
\end{equation*}
где
$$
\mn = \mn(\alpha -\lambda -(p^{-1} -q^{-1})_+ \e),
\cmn = \cmn(\alpha -\lambda -(p^{-1} -q^{-1})_+ \e).
$$

Доказательство.

Полагая $ \mathcal J = \{j \in \Nu_{1,d}^1:
\alpha_j -\lambda_j -(p^{-1} -q^{-1})_+ = \mn \}, $
фиксируем вектор $ \beta \in \R_+^d, $ обладающий следующими свойствами:
$ \beta_j =1, j \in \mathcal J, \beta_j >1 $ и $ \beta_j^{-1}
( \alpha_j -\lambda_j -(p^{-1} -q^{-1})_+) > \mn,
j \in \Nu_{1,d}^1 \setminus \mathcal J. $

Далее, при $ l = l(\alpha) $ определим семейство точек
\begin{multline*}
\mathtt x_{\kappa, \nu}^{d, l -\e, \rho} = \xi_{\delta,
x^0}^{d,l-\e,\rho} = x^0 +\delta \xi_{\e,0}^{d,l-\e,\rho}
\text{(см. п. 2.1),  при } x^0 = 2^{-\kappa} \nu, \delta \\
= 2^{-\kappa},
\rho \in \Z_+^d(l -\e), \nu \in \Nu_{0, 2^\kappa -\e}^d, \kappa \in \Z_+^d,
\end{multline*}
и для $ r \in \N $ рассмотрим множество точек
$$
\{ \mathtt x_{\kappa, \nu}^{d, l -\e, \rho}:
\rho \in \Z_+^d(l -\e), \nu \in \Nu_{0, 2^\kappa -\e}^d,
\kappa \in \Z_+^d, (\kappa, \beta) \le r \}.
$$
Число этих точек, учитывая (1.2.1), удовлетворяет неравенству
\begin{multline*}
\card \{ \mathtt x_{\kappa, \nu}^{d, l -\e, \rho}: \rho \in
\Z_+^d(l -\e), \nu \in \Nu_{0, 2^\kappa -\e}^d, \kappa \in \Z_+^d,
(\kappa, \beta) \le r \}\\
 \le \sum_{\kappa \in \Z_+^d: (\kappa,
\beta) \le r} l^{\e} 2^{(\kappa, \e)} \le l^{\e} c_{28}
2^{\mx(\beta^{-1} \e) r} r^{\cmx(\beta^{-1} \e) -1} = l^{\e}
c_{28} 2^r r^{\cmn -1} = c_{29} 2^r r^{\cmn -1}.
\end{multline*}

Для $ n \ge n_0 = 2 c_{29} $ выберем
$ r \in \N $ так, чтобы соблюдалось соотношение
\begin{equation*} \tag{2.2.21}
c_{29} 2^r r^{\cmn -1} \le n < c_{29} 2^{r+1} (r+1)^{\cmn -1},
\end{equation*}
и построим соответствующую этому $ r \in \N $ систему точек
$$
\{ \mathtt x_{\kappa, \nu}^{d, l -\e, \rho}: \rho \in
\Z_+^d(l -\e), \nu \in \Nu_{0, 2^\kappa -\e}^d, \kappa \in \Z_+^d,
(\kappa, \beta) \le r \}.
$$
Тогда в условиях теоремы, взяв некоторое $ m \in \Z_+^d $ так, чтобы $ \lambda
\in \Z_+^d(m), $ нетрудно видеть, что существуют отображения
$ \phi \in \Phi_n(C(I^d)) $ и $ A \in \overline {\mathcal A}^n(L_q(I^d)) $
такие, что для $ f \in \mathcal S_{p,\theta}^\alpha \mathcal B(I^d) $
имеет место представление
$$
A \circ \phi (f) = \sum_{ \kappa \in \Z_+^d: (\kappa, \beta) \le r}
\D^\lambda \mathcal R_\kappa^{d, l -\e,m}f.
$$
При этом, в силу (2.2.18), (2.2.1) для
$ f \in \mathcal S_{p,\theta}^\alpha \mathcal B(I^d) $ выполняется неравенство
\begin{multline*} \tag{2.2.22}
\| \D^\lambda f -A \circ \phi(f) \|_{L_q(I^d)} = \| \D^\lambda f -
\sum_{ \kappa \in \Z_+^d: (\kappa, \beta) \le r} \D^\lambda
\mathcal R_\kappa^{d, l -\e,m}f \|_{L_q(I^d)} \\
= \| \sum_{ \kappa
\in \Z_+^d: (\kappa, \beta) > r} \D^\lambda \mathcal R_\kappa^{d,
l -\e,m}f \|_{L_q(I^d)} \le \sum_{ \kappa \in \Z_+^d: (\kappa,
\beta) > r} \| \D^\lambda \mathcal R_\kappa^{d, l -\e,m}f
\|_{L_q(I^d)} \\
\le \sum_{ \kappa \in \Z_+^d: (\kappa, \beta) > r}
c_{1} 2^{(\kappa, \lambda +(p^{-1} -q^{-1})_+ \e)} \biggl(\Omega^{l
\chi_{\s(\kappa)}}(f, (c_{2} 2^{-\kappa})^{\s(\kappa)})_{L_p(I^d)}\\
+\sum_{J \subset \Nu_{1,d}^1: J \ne \emptyset} \biggl(\int_{ (c_{3}
I^d)^J} (\prod_{j \in J} u_j^{-p(\alpha_j -\epsilon_j) -1})
(\Omega^{l \chi_{J \cup \s(\kappa)}}( f, (u \chi_{J \setminus
\s(\kappa)} \\
+ 2^{-\kappa} u \chi_{\s(\kappa) \cap J} +c_{2}
2^{-\kappa} \chi_{\s(\kappa) \setminus J})^{J \cup
\s(\kappa)})_{L_p(I^d)})^p du^J\biggr)^{1/p}\biggr) \\
= \sum_{ \kappa \in
\Z_+^d: (\kappa, \beta) > r} c_{1} 2^{(\kappa, \lambda +(p^{-1}
-q^{-1})_+ \e)} \Omega^{l \chi_{\s(\kappa)}}(f, (c_{2}
2^{-\kappa})^{\s(\kappa)})_{L_p(I^d)} \\
+\sum_{ \kappa \in \Z_+^d:
(\kappa, \beta) > r} \sum_{J \subset \Nu_{1,d}^1: J \ne \emptyset}
c_{1} 2^{(\kappa, \lambda +(p^{-1} -q^{-1})_+ \e)} \biggl(\int_{ (c_{3}
I^d)^J} (\prod_{j \in J} u_j^{-p(\alpha_j -\epsilon_j) -1})\\
\times(\Omega^{l \chi_{J \cup \s(\kappa)}}( f, (u \chi_{J \setminus
\s(\kappa)} + 2^{-\kappa} u \chi_{\s(\kappa) \cap J} +c_{2}
2^{-\kappa} \chi_{\s(\kappa) \setminus J})^{J \cup
\s(\kappa)})_{L_p(I^d)})^p du^J\biggr)^{1/p} \\
= c_{1} \biggl(\sum_{ \kappa \in
\Z_+^d: (\kappa, \beta) > r} 2^{(\kappa, \lambda +(p^{-1}
-q^{-1})_+ \e)} \Omega^{l \chi_{\s(\kappa)}}(f, (c_{2}
2^{-\kappa})^{\s(\kappa)})_{L_p(I^d)} \\
+\sum_{J \subset
\Nu_{1,d}^1: J \ne \emptyset} \sum_{ \kappa \in \Z_+^d: (\kappa,
\beta) > r} 2^{(\kappa, \lambda +(p^{-1} -q^{-1})_+ \e)} \biggl(\int_{
(c_{3} I^d)^J} (\prod_{j \in J} u_j^{-p(\alpha_j -\epsilon_j) -1})\\
(\Omega^{l \chi_{J \cup \s(\kappa)}}( f, (u \chi_{J \setminus
\s(\kappa)} + 2^{-\kappa} u \chi_{\s(\kappa) \cap J} +c_{2}
2^{-\kappa} \chi_{\s(\kappa) \setminus J})^{J \cup
\s(\kappa)})_{L_p(I^d)})^p du^J\biggr)^{1/p}\biggr)\\
 = c_{1} \biggl(\sum_{ \kappa \in
\Z_+^d: (\kappa, \beta) > r} 2^{(\kappa, \lambda +(p^{-1}
-q^{-1})_+ \e)} \Omega^{l \chi_{\s(\kappa)}}(f, (c_{2}
2^{-\kappa})^{\s(\kappa)})_{L_p(I^d)}\\
+\sum_{J \subset
\Nu_{1,d}^1: J \ne \emptyset} \sum_{J^\prime \subset \Nu_{1,d}^1:
J^\prime \ne \emptyset} \sum_{ \kappa \in \Z_+^d: (\kappa, \beta)
> r, \s(\kappa) = J^\prime} 2^{(\kappa, \lambda +(p^{-1}
-q^{-1})_+ \e)} \\
\times\biggl(\int_{ (c_{3} I^d)^J} (\prod_{j \in J}
u_j^{-p(\alpha_j -\epsilon_j) -1}) (\Omega^{l \chi_{J \cup
J^\prime}}( f, (u \chi_{J \setminus J^\prime} \\
+ 2^{-\kappa} u
\chi_{J^\prime \cap J} +c_{2} 2^{-\kappa} \chi_{J^\prime \setminus
J})^{J \cup J^\prime})_{L_p(I^d)})^p du^J\biggr)^{1/p}\biggr).
\end{multline*}

Оценивая первую сумму в правой части (2.2.22), с помощью неравенства
Г\"ельдера для $ f \in K, r \in \N $ получаем
\begin{multline*} \tag{2.2.23}
\sum_{ \kappa \in \Z_+^d: (\kappa, \beta) > r}
2^{(\kappa, \lambda +(p^{-1} -q^{-1})_+ \e)}
\Omega^{l \chi_{\s(\kappa)}}(f, (c_{2} 2^{-\kappa})^{\s(\kappa)})_{L_p(I^d)}\\
=\sum_{ \kappa \in \Z_+^d: (\kappa, \beta) > r}
2^{-(\kappa, \alpha -\lambda -(p^{-1} -q^{-1})_+ \e)}
2^{(\kappa, \alpha)}
\Omega^{l \chi_{\s(\kappa)}}(f, (c_{2} 2^{-\kappa})^{\s(\kappa)})_{L_p(I^d)}\\
\le \biggl(\sum_{ \kappa \in \Z_+^d: (\kappa, \beta) > r}
2^{-(\kappa, \alpha -\lambda -(p^{-1}  -q^{-1})_+ \e)
\theta^\prime}\biggr)^{1/\theta^\prime}\\
\times\biggl(\sum_{ \kappa \in \Z_+^d: (\kappa, \beta) > r}
(2^{(\kappa, \alpha)}
\Omega^{l \chi_{\s(\kappa)}}(f,
(c_{2} 2^{-\kappa})^{\s(\kappa)})_{L_p(I^d)})^\theta \biggr)^{1/\theta},
\theta^\prime = \theta/(\theta -1).
\end{multline*}

Пользуясь (1.2.2), выводим
\begin{multline*} \tag{2.2.24}
\biggl(\sum_{ \kappa \in \Z_+^d: (\kappa, \beta) > r}
2^{-(\kappa, \alpha -\lambda -(p^{-1}  -q^{-1})_+ \e)
\theta^\prime}\biggr)^{1/\theta^\prime} \\
\le
(c_{30} 2^{-\mn(\theta^\prime \beta^{-1} (\alpha -\lambda -(p^{-1}
-q^{-1})_+ \e)) r} r^{\cmn(\theta^\prime
\beta^{-1} (\alpha -\lambda -(p^{-1}
-q^{-1})_+ \e)) -1})^{1/\theta^\prime} \\
=
(c_{30} 2^{-\theta^\prime \mn(\alpha -\lambda -(p^{-1}
-q^{-1})_+ \e) r} r^{\cmn(\alpha -\lambda -(p^{-1}
-q^{-1})_+ \e) -1})^{1/\theta^\prime} \\
=
c_{31} 2^{-\mn(\alpha -\lambda -(p^{-1}  -q^{-1})_+ \e) r}
r^{(1 -1/\theta)(\cmn(\alpha -\lambda -(p^{-1} -q^{-1})_+ \e) -1)}, r \in \N.
\end{multline*}

Как показано в [5], существует константа $ c_{32}(d,\alpha,p,\theta) >0 $
такая, что для $ f \in \mathcal S_{p,\theta}^\alpha \mathcal B(I^d) $
при $ r \in \N $ имеет место неравенство
\begin{equation*} \tag{2.2.25}
\biggl( \sum_{\kappa \in \Z_+^d: (\kappa, \beta) > r} (2^{(\kappa,
\alpha)} \Omega^{l \chi_{\s(\kappa)}}(f, c_{2}
(2^{-\kappa})^{\s(\kappa)})_{L_p(I^d)})^\theta \biggr)^{1/\theta} \le
c_{32}.
\end{equation*}

Подставляя (2.2.24) и (2.2.25) в (2.2.23), приходим к неравенству
\begin{multline*} \tag{2.2.26}
\sum_{ \kappa \in \Z_+^d: (\kappa, \beta) > r}
2^{(\kappa, \lambda +(p^{-1} -q^{-1})_+ \e)}
\Omega^{l \chi_{\s(\kappa)}}(f, (c_{2} 2^{-\kappa})^{\s(\kappa)})_{L_p(I^d)}\\
\le
c_{33} 2^{-\mn(\alpha -\lambda -(p^{-1}  -q^{-1})_+ \e) r}
r^{(1 -1/\theta)(\cmn(\alpha -\lambda -(p^{-1} -q^{-1})_+ \e) -1)},
f \in K, r \in \N.
\end{multline*}

Продолжая оценку правой части (2.2.22), рассмотрим два случая.
В первом случае, когда $ \theta > p, $ для
$ f \in K, J \subset \Nu_{1,d}^1: J \ne \emptyset,
 J^\prime \subset \Nu_{1,d}^1: J^\prime \ne \emptyset, r \in \N, $ благодаря
неравенству Г\"ельдера, с уч\"етом (2.2.24) имеем
\begin{multline*} \tag{2.2.27}
\sum_{ \kappa \in \Z_+^d: (\kappa, \beta) > r, \s(\kappa) = J^\prime}
2^{(\kappa, \lambda +(p^{-1} -q^{-1})_+ \e)}
\biggl(\int_{ (c_{3} I^d)^J}
(\prod_{j \in J}
u_j^{-p(\alpha_j -\epsilon_j) -1})\\
\times(\Omega^{l \chi_{J \cup J^\prime}}( f, (u \chi_{J \setminus J^\prime} +
2^{-\kappa} u \chi_{J^\prime \cap J} +c_{2} 2^{-\kappa}
\chi_{J^\prime \setminus J})^{J \cup J^\prime})_{L_p(I^d)})^p du^J\biggr)^{1/p} \\
=
\sum_{ \kappa \in \Z_+^d: (\kappa, \beta) > r, \s(\kappa) = J^\prime}
2^{-(\kappa, \alpha -\lambda -(p^{-1}  -q^{-1})_+ \e)}
\biggl(\int_{ (c_{3} I^d)^J}
(\prod_{j \in J}
u_j^{-p(\alpha_j -\epsilon_j) -1})\\
\times( 2^{(\kappa^{J^\prime}, \alpha^{J^\prime})}
\Omega^{l \chi_{J \cup J^\prime}}( f, (u \chi_{J \setminus J^\prime} +
2^{-\kappa} u \chi_{J^\prime \cap J} +c_{2} 2^{-\kappa}
\chi_{J^\prime \setminus J})^{J \cup J^\prime})_{L_p(I^d)})^p du^J\biggr)^{1/p}\\
 \le \biggl(\sum_{ \kappa \in \Z_+^d: (\kappa, \beta) > r, \s(\kappa) = J^\prime}
2^{-(\kappa, \alpha -\lambda -(p^{-1}  -q^{-1})_+ \e)
\theta^\prime}\biggr)^{1/\theta^\prime}\\
\times\biggl(\sum_{ \kappa \in \Z_+^d: (\kappa, \beta) > r, \s(\kappa) = J^\prime}
(\int_{ (c_{3} I^d)^J}
(\prod_{j \in J}
u_j^{-p(\alpha_j -\epsilon_j) -1})\\
\times( 2^{(\kappa^{J^\prime}, \alpha^{J^\prime})}
\Omega^{l \chi_{J \cup J^\prime}}( f, (u \chi_{J \setminus J^\prime} +
2^{-\kappa} u \chi_{J^\prime \cap J} +c_{2} 2^{-\kappa}
\chi_{J^\prime \setminus J})^{J \cup J^\prime})_{L_p(I^d)})^p
du^J)^{\theta/p}\biggr)^{1/\theta}\\
\le
\biggl(\sum_{ \kappa \in \Z_+^d: (\kappa, \beta) > r}
2^{-(\kappa, \alpha -\lambda -(p^{-1}  -q^{-1})_+ \e)
\theta^\prime}\biggr)^{1/\theta^\prime}\\
\biggl(\sum_{ \kappa \in \Z_+^d: \s(\kappa) = J^\prime}
(\int_{ (c_{3} I^d)^J}
(\prod_{j \in J}
u_j^{-p(\alpha_j -\epsilon_j) -1})\\
\times( 2^{(\kappa^{J^\prime}, \alpha^{J^\prime})}
\Omega^{l \chi_{J \cup J^\prime}}( f, (u \chi_{J \setminus J^\prime} +
2^{-\kappa} u \chi_{J^\prime \cap J} +c_{2} 2^{-\kappa}
\chi_{J^\prime \setminus J})^{J \cup J^\prime})_{L_p(I^d)})^p
du^J)^{\theta/p} \biggr)^{1/\theta}\\
\le c_{34} 2^{-\mn(\alpha -\lambda -(p^{-1}  -q^{-1})_+ \e) r}
r^{(1 -1/\theta)(\cmn(\alpha -\lambda -(p^{-1} -q^{-1})_+ \e) -1)}\\
\times\biggl(\sum_{ \kappa \in \Z_+^d: \s(\kappa) = J^\prime}
\biggl(\int_{ (c_{3} I^d)^J}
(\prod_{j \in J}
u_j^{-p(\alpha_j -\epsilon_j) -1})
( 2^{(\kappa^{J^\prime}, \alpha^{J^\prime})}
\Omega^{l \chi_{J \cup J^\prime}}( f, (u \chi_{J \setminus J^\prime} \\
+2^{-\kappa} u \chi_{J^\prime \cap J} +c_{2} 2^{-\kappa}
\chi_{J^\prime \setminus J})^{J \cup J^\prime})_{L_p(I^d)})^p
du^J\biggr)^{\theta/p} \biggr)^{1/\theta}.
\end{multline*}

Оценивая слагаемые в правой части (2.2.27), пользуясь неравенством
Г\"ельдера с показателем $ \theta/p >1, $ для $ f \in K $
при $ J \subset \Nu_{1,d}^1: J \ne \emptyset,
J^\prime \subset \Nu_{1,d}^1: J^\prime \ne \emptyset,
\kappa \in \Z_+^d: \s(\kappa) = J^\prime, $ выводим
\begin{multline*} \tag{2.2.28}
\int_{ (c_{3} I^d)^J}
(\prod_{j \in J}
u_j^{-p(\alpha_j -\epsilon_j) -1})
( 2^{(\kappa^{J^\prime}, \alpha^{J^\prime})}\\
\times\Omega^{l \chi_{J \cup J^\prime}}( f, (u \chi_{J \setminus J^\prime} +
2^{-\kappa} u \chi_{J^\prime \cap J} +c_{2} 2^{-\kappa}
\chi_{J^\prime \setminus J})^{J \cup J^\prime})_{L_p(I^d)})^p du^J \\
=\int_{ (c_{3} I^d)^J}
(\prod_{j \in J}u_j^{\frac{1}{2} p \epsilon_j -(\theta -p)/\theta})
(\prod_{j \in J}u_j^{-p(\alpha_j -\frac{1}{2} \epsilon_j) -p/\theta})
( 2^{(\kappa^{J^\prime}, \alpha^{J^\prime})}\\
\times\Omega^{l \chi_{J \cup J^\prime}}( f, (u \chi_{J \setminus J^\prime} +
2^{-\kappa} u \chi_{J^\prime \cap J} +c_{2} 2^{-\kappa}
\chi_{J^\prime \setminus J})^{J \cup J^\prime})_{L_p(I^d)})^p du^J\\
\le \biggl(\int_{ (c_{3} I^d)^J}(\prod_{j \in J}
u_j^{\frac{1}{2} p \epsilon_j -(\theta -p)/\theta})^{\theta/(\theta -p)}
du^J\biggr)^{(\theta -p)/\theta}\\
\times\biggl(\int_{ (c_{3} I^d)^J}\biggl((\prod_{j \in J}
u_j^{-p(\alpha_j -\frac{1}{2} \epsilon_j) -p/\theta})
( 2^{(\kappa^{J^\prime}, \alpha^{J^\prime})}\\
\times\Omega^{l \chi_{J \cup J^\prime}}( f, (u \chi_{J \setminus J^\prime} +
2^{-\kappa} u \chi_{J^\prime \cap J} +c_{2} 2^{-\kappa}
\chi_{J^\prime \setminus J})^{J \cup J^\prime})_{L_p(I^d)})^p \biggr)^{\theta /p}
du^J\biggr)^{p/\theta}\\
=\biggl(\int_{ (c_{3} I^d)^J}(\prod_{j \in J}
u_j^{\frac{1}{2} p \epsilon_j \theta/(\theta -p) -1})
du^J\biggr)^{(\theta -p)/\theta}
\biggl(\int_{ (c_{3} I^d)^J}(\prod_{j \in J}
u_j^{-\theta (\alpha_j -\frac{1}{2} \epsilon_j) -1})
( 2^{(\kappa^{J^\prime}, \alpha^{J^\prime})}\\
\times\Omega^{l \chi_{J \cup J^\prime}}( f, (u \chi_{J \setminus J^\prime} +
2^{-\kappa} u \chi_{J^\prime \cap J} +c_{2} 2^{-\kappa}
\chi_{J^\prime \setminus J})^{J \cup J^\prime})_{L_p(I^d)} )^\theta
du^J\biggr)^{p/\theta}\\
=c_{35}\biggl(\int_{ (c_{3} I^d)^J}
(\prod_{j \in J}
u_j^{-\theta (\alpha_j -\frac{1}{2} \epsilon_j) -1})
( 2^{(\kappa^{J^\prime}, \alpha^{J^\prime})}\\
\times\Omega^{l \chi_{J \cup J^\prime}}( f, (u \chi_{J \setminus J^\prime} +
2^{-\kappa} u \chi_{J^\prime \cap J} +c_{2} 2^{-\kappa}
\chi_{J^\prime \setminus J})^{J \cup J^\prime})_{L_p(I^d)} )^\theta
du^J\biggr)^{p/\theta}.
\end{multline*}

Из (2.2.28) с учетом замечания после леммы 1.2.1 и теоремы Лебега о
предельном переходе под знаком интеграла вытекает, что
при $ J \subset \Nu_{1,d}^1: J \ne \emptyset,
J^\prime \subset \Nu_{1,d}^1: J^\prime \ne \emptyset, $ для $ f \in K $
выполняется неравенство
\begin{multline*} \tag{2.2.29}
\biggl(\sum_{ \kappa \in \Z_+^d: \s(\kappa) = J^\prime}
\biggl(\int_{ (c_{3} I^d)^J}
(\prod_{j \in J}u_j^{-p(\alpha_j -\epsilon_j) -1})
( 2^{(\kappa^{J^\prime}, \alpha^{J^\prime})}\\
\times\Omega^{l \chi_{J \cup J^\prime}}( f, (u \chi_{J \setminus J^\prime} +
2^{-\kappa} u \chi_{J^\prime \cap J} +c_{2} 2^{-\kappa}
\chi_{J^\prime \setminus J})^{J \cup J^\prime})_{L_p(I^d)})^p
du^J\biggr)^{\theta/p} \biggr)^{1/\theta}\\
\le c_{36}\biggl(\sum_{ \kappa \in \Z_+^d: \s(\kappa) = J^\prime}
\int_{ (c_{3} I^d)^J}
(\prod_{j \in J}
u_j^{-\theta (\alpha_j -\frac{1}{2} \epsilon_j) -1})
( 2^{(\kappa^{J^\prime}, \alpha^{J^\prime})}\\
\times\Omega^{l \chi_{J \cup J^\prime}}( f, (u \chi_{J \setminus J^\prime} +
2^{-\kappa} u \chi_{J^\prime \cap J} +c_{2} 2^{-\kappa}
\chi_{J^\prime \setminus J})^{J \cup J^\prime})_{L_p(I^d)} )^\theta
du^J \biggr)^{1/\theta}\\
=c_{36}\biggl(\int_{ (c_{3} I^d)^J}
(\prod_{j \in J}
u_j^{-\theta (\alpha_j -\frac{1}{2} \epsilon_j) -1})
\biggl(\sum_{ \kappa \in \Z_+^d: \s(\kappa) = J^\prime}
( 2^{(\kappa^{J^\prime}, \alpha^{J^\prime})}\\
\times\Omega^{l \chi_{J \cup J^\prime}}( f, (u \chi_{J \setminus J^\prime} +
2^{-\kappa} u \chi_{J^\prime \cap J} +c_{2} 2^{-\kappa}
\chi_{J^\prime \setminus J})^{J \cup J^\prime})_{L_p(I^d)} )^\theta\biggr)
du^J \biggr)^{1/\theta}.
\end{multline*}

Далее, учитывая, что при  $ J \subset \Nu_{1,d}^1: J \ne \emptyset,
J^\prime \subset \Nu_{1,d}^1: J^\prime \ne \emptyset,
\kappa \in \Z_+^d: \s(\kappa) = J^\prime,
u^J \in (\R_+^d)^J $ для $ f \in K $ соблюдается неравенство
\begin{multline*}
( 2^{(\kappa^{J^\prime}, \alpha^{J^\prime})}
\Omega^{l \chi_{J \cup J^\prime}}( f, (u \chi_{J \setminus J^\prime} +
2^{-\kappa} u \chi_{J^\prime \cap J} +c_{2} 2^{-\kappa}
\chi_{J^\prime \setminus J})^{J \cup J^\prime})_{L_p(I^d)} )^\theta \\
=(u^{J \setminus J^\prime})^{\theta \alpha^{J \setminus J^\prime}}
((u^{J \setminus J^\prime})^{-\alpha^{J \setminus J^\prime}}
(2^{-\kappa^{J^\prime}})^{-\alpha^{J^\prime}}\\
\times\Omega^{l \chi_{J \cup J^\prime}}( f, (u \chi_{J \setminus J^\prime} +
2^{-\kappa} u \chi_{J^\prime \cap J} +c_{2} 2^{-\kappa}
\chi_{J^\prime \setminus J})^{J \cup J^\prime})_{L_p(I^d)} )^\theta \\
=(u^{J \setminus J^\prime})^{\theta \alpha^{J \setminus J^\prime}}
\int_{(u +u I^d)^{J \setminus J^\prime} \times (2^{-\kappa} +
2^{-\kappa} I^d)^{J^\prime}}
(u^{J \setminus J^\prime})^{-\e^{J \setminus J^\prime} -\theta
\alpha^{J \setminus J^\prime}}
(2^{-\kappa^{J^\prime}})^{-\e^{J^\prime} -\theta \alpha^{J^\prime}}\\
\times(\Omega^{l \chi_{J \cup J^\prime}}( f, (u \chi_{J \setminus J^\prime} +
2^{-\kappa} u \chi_{J^\prime \cap J} +c_{2} 2^{-\kappa}
\chi_{J^\prime \setminus J})^{J \cup J^\prime})_{L_p(I^d)} )^\theta
d\tau^{J \cup J^\prime} \\
\le c_{37}(u^{J \setminus J^\prime})^{\theta \alpha^{J \setminus J^\prime}}
\int_{(u +u I^d)^{J \setminus J^\prime} \times (2^{-\kappa} +
2^{-\kappa} I^d)^{J^\prime}}
(\tau^{J \setminus J^\prime})^{-\e^{J \setminus J^\prime} -\theta
\alpha^{J \setminus J^\prime}}
(\tau^{J^\prime})^{-\e^{J^\prime} -\theta \alpha^{J^\prime}}\\
\times(\Omega^{l \chi_{J \cup J^\prime}}( f, (\tau \chi_{J \setminus J^\prime} +
\tau u \chi_{J^\prime \cap J} +c_{2} \tau
\chi_{J^\prime \setminus J})^{J \cup J^\prime})_{L_p(I^d)} )^\theta
d\tau^{J \cup J^\prime}\\
 \le c_{37}
(u^{J \setminus J^\prime})^{\theta \alpha^{J \setminus J^\prime}}
\int_{(\R_+^d)^{J \setminus J^\prime} \times (2^{-\kappa} +
2^{-\kappa} I^d)^{J^\prime}}
(\tau^{J \cup J^\prime})^{-\e^{J \cup J^\prime} -\theta
\alpha^{J \cup J^\prime}}\\
\times(\Omega^{l \chi_{J \cup J^\prime}}( f, (\tau \chi_{J \setminus J^\prime} +
\tau u \chi_{J^\prime \cap J} +c_{2} \tau
\chi_{J^\prime \setminus J})^{J \cup J^\prime})_{L_p(I^d)} )^\theta
d\tau^{J \cup J^\prime},
\end{multline*}
заключаем, что
\begin{multline*} \tag{2.2.30}
\sum_{ \kappa \in \Z_+^d: \s(\kappa) = J^\prime}
( 2^{(\kappa^{J^\prime}, \alpha^{J^\prime})}
\Omega^{l \chi_{J \cup J^\prime}}( f, (u \chi_{J \setminus J^\prime} +
2^{-\kappa} u \chi_{J^\prime \cap J} +c_{2} 2^{-\kappa}
\chi_{J^\prime \setminus J})^{J \cup J^\prime})_{L_p(I^d)} )^\theta \\
\le c_{37}(u^{J \setminus J^\prime})^{\theta \alpha^{J \setminus J^\prime}}
\sum_{ \kappa \in \Z_+^d: \s(\kappa) = J^\prime}
\int_{(\R_+^d)^{J \setminus J^\prime} \times (2^{-\kappa} +
2^{-\kappa} I^d)^{J^\prime}}
(\tau^{J \cup J^\prime})^{-\e^{J \cup J^\prime} -\theta
\alpha^{J \cup J^\prime}}\\
\times(\Omega^{l \chi_{J \cup J^\prime}}( f, (\tau \chi_{J \setminus J^\prime} +
\tau u \chi_{J^\prime \cap J} +c_{2} \tau
\chi_{J^\prime \setminus J})^{J \cup J^\prime})_{L_p(I^d)} )^\theta
d\tau^{J \cup J^\prime} \\
=c_{37}(u^{J \setminus J^\prime})^{\theta \alpha^{J \setminus J^\prime}}
\int_{\cup_{ \kappa \in \Z_+^d: \s(\kappa) = J^\prime}
(\R_+^d)^{J \setminus J^\prime} \times (2^{-\kappa} +
2^{-\kappa} I^d)^{J^\prime}}
(\tau^{J \cup J^\prime})^{-\e^{J \cup J^\prime} -\theta
\alpha^{J \cup J^\prime}}\\
\times(\Omega^{l \chi_{J \cup J^\prime}}( f, (\tau \chi_{J \setminus J^\prime} +
\tau u \chi_{J^\prime \cap J} +c_{2} \tau
\chi_{J^\prime \setminus J})^{J \cup J^\prime})_{L_p(I^d)} )^\theta
d\tau^{J \cup J^\prime}\\
 =c_{37}
(u^{J \setminus J^\prime})^{\theta \alpha^{J \setminus J^\prime}}
\int_{ (\R_+^d)^{J \setminus J^\prime} \times (I^d)^{J^\prime}}
(\tau^{J \cup J^\prime})^{-\e^{J \cup J^\prime} -\theta
\alpha^{J \cup J^\prime}}\\
\times(\Omega^{l \chi_{J \cup J^\prime}}( f, (\tau \chi_{J \setminus J^\prime} +
\tau u \chi_{J^\prime \cap J} +c_{2} \tau
\chi_{J^\prime \setminus J})^{J \cup J^\prime})_{L_p(I^d)} )^\theta
d\tau^{J \cup J^\prime} \\
\le c_{37}
(u^{J \setminus J^\prime})^{\theta \alpha^{J \setminus J^\prime}}
\int_{ (\R_+^d)^{J \cup J^\prime}}
(\tau^{J \cup J^\prime})^{-\e^{J \cup J^\prime} -\theta
\alpha^{J \cup J^\prime}}\\
\times(\Omega^{l \chi_{J \cup J^\prime}}( f, (\tau \chi_{J \setminus J^\prime} +
\tau u \chi_{J^\prime \cap J} +c_{2} \tau
\chi_{J^\prime \setminus J})^{J \cup J^\prime})_{L_p(I^d)} )^\theta
d\tau^{J \cup J^\prime} \\
=c_{37}(u^{J \setminus J^\prime})^{\theta \alpha^{J \setminus J^\prime}}
\int_{ (\R_+^d)^{J \cup J^\prime}}
(\prod_{j \in (J \setminus J^\prime)}
t_j^{-1 -\theta \alpha_j})
(\prod_{j \in (J \cap J^\prime)}
(t_j /u_j)^{-1 -\theta \alpha_j})\\
\times(\prod_{j \in (J^\prime \setminus J)}
(t_j /c_{2})^{-1 -\theta \alpha_j})
(\Omega^{l \chi_{J \cup J^\prime}}( f,
t^{J \cup J^\prime})_{L_p(I^d)} )^\theta
(\prod_{j \in (J \cap J^\prime)}
(1 /u_j))
(\prod_{j \in (J^\prime \setminus J)}
(1 /c_{2}))
dt^{J \cup J^\prime} \\
=c_{38}
(u^{J \setminus J^\prime})^{\theta \alpha^{J \setminus J^\prime}}
\int_{ (\R_+^d)^{J \cup J^\prime}}
(\prod_{j \in (J \cap J^\prime)}
u_j^{\theta \alpha_j})
(t^{J \cup J^\prime})^{-\e^{J \cup J^\prime} -\theta
\alpha^{J \cup J^\prime}}\\
\times(\Omega^{l \chi_{J \cup J^\prime}}( f,
t^{J \cup J^\prime})_{L_p(I^d)} )^\theta
dt^{J \cup J^\prime}
=c_{38}
(u^J)^{\theta \alpha^J}
\int_{ (\R_+^d)^{J \cup J^\prime}}
(t^{J \cup J^\prime})^{-\e^{J \cup J^\prime} -\theta
\alpha^{J \cup J^\prime}}\\
\times(\Omega^{l \chi_{J \cup J^\prime}}( f,
t^{J \cup J^\prime})_{L_p(I^d)} )^\theta
dt^{J \cup J^\prime}
\le c_{38} (u^J)^{\theta \alpha^J}.
\end{multline*}

Подставляя оценку (2.2.30) в (2.2.29), находим, что
\begin{multline*} \tag{2.2.31}
\biggl(\sum_{ \kappa \in \Z_+^d: \s(\kappa) = J^\prime}
\biggl(\int_{ (c_{3} I^d)^J}
(\prod_{j \in J}
u_j^{-p(\alpha_j -\epsilon_j) -1})
( 2^{(\kappa^{J^\prime}, \alpha^{J^\prime})}\\
\times\Omega^{l \chi_{J \cup J^\prime}}( f, (u \chi_{J \setminus J^\prime} +
2^{-\kappa} u \chi_{J^\prime \cap J} +c_{2} 2^{-\kappa}
\chi_{J^\prime \setminus J})^{J \cup J^\prime})_{L_p(I^d)})^p
du^J\biggr)^{\theta/p} \biggr)^{1/\theta}\\
\le c_{36}\biggl(\int_{ (c_{3} I^d)^J}
(\prod_{j \in J}u_j^{-\theta (\alpha_j -\frac{1}{2} \epsilon_j) -1})
c_{38} (u^J)^{\theta \alpha^J}
du^J\biggr)^{1/\theta} \\
=c_{39}\biggl(\int_{ (c_{3} I^d)^J}
(u^J)^{\frac{1}{2} \theta \epsilon^J -\e^J}
du^J\biggr)^{1/\theta} = c_{40}, \\
J \subset \Nu_{1,d}^1: J \ne \emptyset,
J^\prime \subset \Nu_{1,d}^1: J^\prime \ne \emptyset, f \in K.
\end{multline*}

Соединяя (2.2.27) с (2.2.31), приходим к неравенству
\begin{multline*} \tag{2.2.32}
\sum_{ \kappa \in \Z_+^d: (\kappa, \beta) > r, \s(\kappa) = J^\prime}
2^{(\kappa, \lambda +(p^{-1} -q^{-1})_+ \e)}
\biggl(\int_{ (c_{3} I^d)^J}
(\prod_{j \in J}
u_j^{-p(\alpha_j -\epsilon_j) -1})\\
\times(\Omega^{l \chi_{J \cup J^\prime}}( f, (u \chi_{J \setminus J^\prime} +
2^{-\kappa} u \chi_{J^\prime \cap J} +c_{2} 2^{-\kappa}
\chi_{J^\prime \setminus J})^{J \cup J^\prime})_{L_p(I^d)})^p du^J\biggr)^{1/p} \\
\le
c_{41} 2^{-\mn(\alpha -\lambda -(p^{-1}  -q^{-1})_+ \e) r}
r^{(1 -1/\theta)(\cmn(\alpha -\lambda -(p^{-1} -q^{-1})_+ \e) -1)},\\
f \in K, r \in \N, J \subset \Nu_{1,d}^1: J \ne \emptyset,
J^\prime \subset \Nu_{1,d}^1: J^\prime \ne \emptyset.
\end{multline*}

После подстановки (2.2.26) и (2.2.32) в (2.2.22) получаем, что для
$ f \in K $ и $ r \in \N, $ удовлетворяющего (2.2.21),
при $ \theta > p $ выполняется неравенство
\begin{multline*} \tag{2.2.33}
\| \D^\lambda f -A \circ \phi(f) \|_{L_q(I^d)} \\
\le c_{42} 2^{-\mn(\alpha -\lambda -(p^{-1}  -q^{-1})_+ \e) r}
r^{(1 -1/\theta)(\cmn(\alpha -\lambda -(p^{-1} -q^{-1})_+ \e) -1)}.
\end{multline*}

Проводя оценку правой части (2.2.22) в случае, когда $ \theta \le p, $
для $ f \in K, J \subset \Nu_{1,d}^1: J \ne \emptyset,
J^\prime \subset \Nu_{1,d}^1: J^\prime \ne \emptyset, r \in \N, $ благодаря
неравенству Г\"ельдера, ввиду (2.2.24) с заменой в н\"ем $ \theta $ на $ p, $
а также используя с учетом замечания после леммы 1.2.1 теорему Лебега о
предельном переходе под знаком интеграла, выводим
\begin{multline*} \tag{2.2.34}
\sum_{ \kappa \in \Z_+^d: (\kappa, \beta) > r, \s(\kappa) = J^\prime}
2^{(\kappa, \lambda +(p^{-1} -q^{-1})_+ \e)}
\biggl(\int_{ (c_{3} I^d)^J}
(\prod_{j \in J}
u_j^{-p(\alpha_j -\epsilon_j) -1})\\
\times(\Omega^{l \chi_{J \cup J^\prime}}( f, (u \chi_{J \setminus J^\prime} +
2^{-\kappa} u \chi_{J^\prime \cap J} +c_{2} 2^{-\kappa}
\chi_{J^\prime \setminus J})^{J \cup J^\prime})_{L_p(I^d)})^p du^J\biggr)^{1/p}\\
 =\sum_{ \kappa \in \Z_+^d: (\kappa, \beta) > r, \s(\kappa) = J^\prime}
2^{-(\kappa, \alpha -\lambda -(p^{-1} -q^{-1})_+ \e)}
\biggl(\int_{ (c_{3} I^d)^J}
(\prod_{j \in J}
u_j^{-p(\alpha_j -\epsilon_j) -1})
( 2^{(\kappa^{J^\prime}, \alpha^{J^\prime})}\\
\times\Omega^{l \chi_{J \cup J^\prime}}( f, (u \chi_{J \setminus J^\prime} +
2^{-\kappa} u \chi_{J^\prime \cap J} +c_{2} 2^{-\kappa}
\chi_{J^\prime \setminus J})^{J \cup J^\prime})_{L_p(I^d)})^p du^J\biggr)^{1/p} \\
\le \biggl(\sum_{ \kappa \in \Z_+^d: (\kappa, \beta) > r, \s(\kappa) = J^\prime}
2^{-(\kappa, \alpha -\lambda -(p^{-1} -q^{-1})_+ \e)
p^\prime}\biggr)^{1/p^\prime}\\
\times\biggl(\sum_{ \kappa \in \Z_+^d: (\kappa, \beta) > r, \s(\kappa) = J^\prime}
\int_{ (c_{3} I^d)^J}
(\prod_{j \in J}
u_j^{-p(\alpha_j -\epsilon_j) -1})
( 2^{(\kappa^{J^\prime}, \alpha^{J^\prime})}\\
\times\Omega^{l \chi_{J \cup J^\prime}}( f, (u \chi_{J \setminus J^\prime} +
2^{-\kappa} u \chi_{J^\prime \cap J} +c_{2} 2^{-\kappa}
\chi_{J^\prime \setminus J})^{J \cup J^\prime})_{L_p(I^d)})^p
du^J \biggr)^{1/p}\\
\le
\biggl(\sum_{ \kappa \in \Z_+^d: (\kappa, \beta) > r}
2^{-(\kappa, \alpha -\lambda -(p^{-1} -q^{-1})_+ \e)
p^\prime}\biggr)^{1/p^\prime}\\
\times\biggl(\sum_{ \kappa \in \Z_+^d: \s(\kappa) = J^\prime}
\int_{ (c_{3} I^d)^J}
(\prod_{j \in J}
u_j^{-p(\alpha_j -\epsilon_j) -1})
( 2^{(\kappa^{J^\prime}, \alpha^{J^\prime})}\\
\times\Omega^{l \chi_{J \cup J^\prime}}( f, (u \chi_{J \setminus J^\prime} +
2^{-\kappa} u \chi_{J^\prime \cap J} +c_{2} 2^{-\kappa}
\chi_{J^\prime \setminus J})^{J \cup J^\prime})_{L_p(I^d)})^p
du^J \biggr)^{1/p}\\
\le
c_{43} 2^{-\mn(\alpha -\lambda -(p^{-1}  -q^{-1})_+ \e) r}
r^{(1 -1/p)(\cmn(\alpha -\lambda -(p^{-1} -q^{-1})_+ \e) -1)}\\
\times\biggl(\sum_{ \kappa \in \Z_+^d: \s(\kappa) = J^\prime}
\int_{ (c_{3} I^d)^J}
(\prod_{j \in J}
u_j^{-p(\alpha_j -\epsilon_j) -1})
( 2^{(\kappa^{J^\prime}, \alpha^{J^\prime})}\\
\times\Omega^{l \chi_{J \cup J^\prime}}( f, (u \chi_{J \setminus J^\prime} +
2^{-\kappa} u \chi_{J^\prime \cap J} +c_{2} 2^{-\kappa}
\chi_{J^\prime \setminus J})^{J \cup J^\prime})_{L_p(I^d)})^p
du^J \biggr)^{1/p}\\
=
c_{43} 2^{-\mn(\alpha -\lambda -(p^{-1}  -q^{-1})_+ \e) r}
r^{(1 -1/p)(\cmn(\alpha -\lambda -(p^{-1} -q^{-1})_+ \e) -1)}\\
\times\biggl(\int_{ (c_{3} I^d)^J}
(\prod_{j \in J}
u_j^{-p(\alpha_j -\epsilon_j) -1})
(\sum_{ \kappa \in \Z_+^d: \s(\kappa) = J^\prime}
( 2^{(\kappa^{J^\prime}, \alpha^{J^\prime})}\\
\times\Omega^{l \chi_{J \cup J^\prime}}( f, (u \chi_{J \setminus J^\prime} +
2^{-\kappa} u \chi_{J^\prime \cap J} +c_{2} 2^{-\kappa}
\chi_{J^\prime \setminus J})^{J \cup J^\prime})_{L_p(I^d)})^p )
du^J \biggr)^{1/p}.
\end{multline*}

Применяя неравенство Г\"ельдера с показателем $ \theta /p \le 1, $
а также неравенство (2.2.30), для $ f \in K,  J \subset \Nu_{1,d}^1:
J \ne \emptyset,
J^\prime \subset \Nu_{1,d}^1: J^\prime \ne \emptyset,
u^J \in (\R_+^d)^J $ имеем
\begin{multline*} \tag{2.2.35}
\sum_{ \kappa \in \Z_+^d: \s(\kappa) = J^\prime}
( 2^{(\kappa^{J^\prime}, \alpha^{J^\prime})}
\Omega^{l \chi_{J \cup J^\prime}}( f, (u \chi_{J \setminus J^\prime} +
2^{-\kappa} u \chi_{J^\prime \cap J} +c_{2} 2^{-\kappa}
\chi_{J^\prime \setminus J})^{J \cup J^\prime})_{L_p(I^d)})^p\\
\le(\sum_{ \kappa \in \Z_+^d: \s(\kappa) = J^\prime}
( 2^{(\kappa^{J^\prime}, \alpha^{J^\prime})}
\Omega^{l \chi_{J \cup J^\prime}}( f, (u \chi_{J \setminus J^\prime} +
2^{-\kappa} u \chi_{J^\prime \cap J} +c_{2} 2^{-\kappa}
\chi_{J^\prime \setminus J})^{J \cup J^\prime})_{L_p(I^d)})^\theta )^{p /\theta}\\
\le (c_{38} (u^J)^{\theta \alpha^J} )^{p /\theta} =
c_{44} (u^J)^{p \alpha^J}.
\end{multline*}

Из (2.2.35) вытекает, что при $ J \subset \Nu_{1,d}^1: J \ne \emptyset,
J^\prime \subset \Nu_{1,d}^1: J^\prime \ne \emptyset, \theta \le p, $
для $ f \in K $ имеет место неравенство
\begin{multline*} \tag{2.2.36}
\biggl(\int_{ (c_{3} I^d)^J}
(\prod_{j \in J}
u_j^{-p(\alpha_j -\epsilon_j) -1})
\biggl(\sum_{ \kappa \in \Z_+^d: \s(\kappa) = J^\prime}
( 2^{(\kappa^{J^\prime}, \alpha^{J^\prime})}\\
\times\Omega^{l \chi_{J \cup J^\prime}}( f, (u \chi_{J \setminus J^\prime} +
2^{-\kappa} u \chi_{J^\prime \cap J} +c_{2} 2^{-\kappa}
\chi_{J^\prime \setminus J})^{J \cup J^\prime})_{L_p(I^d)})^p \biggr)
du^J \biggr)^{1/p}\\
\le \biggl(\int_{ (c_{3} I^d)^J}
(\prod_{j \in J}
u_j^{-p(\alpha_j -\epsilon_j) -1})
c_{44} (u^J)^{p \alpha^J} du^J \biggr)^{1/p}\\
 = c_{45} \biggl(\int_{ (c_{3} I^d)^J}
(u^J)^{p \epsilon^J -\e^J} du^J \biggr)^{1/p} = c_{46}.
\end{multline*}

Подставляя оценку (2.2.36) в (2.2.34), находим, что для $ f \in K,
\theta \le p, $ при$ J \subset \Nu_{1,d}^1: J \ne \emptyset,
J^\prime \subset \Nu_{1,d}^1: J^\prime \ne \emptyset, r \in \N $
соблюдается неравенство
\begin{multline*} \tag{2.2.37}
\sum_{ \kappa \in \Z_+^d: (\kappa, \beta) > r, \s(\kappa) = J^\prime}
2^{(\kappa, \lambda +(p^{-1} -q^{-1})_+ \e)}
\biggl(\int_{ (c_{3} I^d)^J}
(\prod_{j \in J}
u_j^{-p(\alpha_j -\epsilon_j) -1})\\
\times(\Omega^{l \chi_{J \cup J^\prime}}( f, (u \chi_{J \setminus J^\prime} +
2^{-\kappa} u \chi_{J^\prime \cap J} +c_{2} 2^{-\kappa}
\chi_{J^\prime \setminus J})^{J \cup J^\prime})_{L_p(I^d)})^p du^J\biggr)^{1/p}\\
\le
c_{47} 2^{-\mn(\alpha -\lambda -(p^{-1}  -q^{-1})_+ \e) r}
r^{(1 -1/p)(\cmn(\alpha -\lambda -(p^{-1} -q^{-1})_+ \e) -1)}.
\end{multline*}

Соединяя (2.2.26) и (2.2.37) с (2.2.22), получаем, что для
$ f \in K $ и $ r \in \N, $ удовлетворяющего (2.2.21),
при $ \theta \le p $ справедливо неравенство
\begin{multline*} \tag{2.2.38}
\| \D^\lambda f -A \circ \phi(f) \|_{L_q(I^d)} \le
c_{48} 2^{-\mn(\alpha -\lambda -(p^{-1}  -q^{-1})_+ \e) r}
r^{(1 -1/\theta)(\cmn(\alpha -\lambda -(p^{-1} -q^{-1})_+ \e) -1)} \\
+c_{49} 2^{-\mn(\alpha -\lambda -(p^{-1}  -q^{-1})_+ \e) r}
r^{(1 -1/p)(\cmn(\alpha -\lambda -(p^{-1} -q^{-1})_+ \e) -1)}\\
\le
c_{50} 2^{-\mn(\alpha -\lambda -(p^{-1}  -q^{-1})_+ \e) r}
r^{(1 -1/p)(\cmn(\alpha -\lambda -(p^{-1} -q^{-1})_+ \e) -1)}.
\end{multline*}

Из (2.2.33), (2.2.38), (2.2.21)   следует (2.2.20). $ \square $
Отметим, что в проведенном доказательстве всюду при выводе
соотношений, содержащих ряды, неявно утверждается, что "если правая часть
выводимого соотношения определена, то левая его часть также определена и
справедливо выводимое соотношение".

\end{document}